\documentclass[11pt,leqno]{article}
\usepackage{amsmath,amssymb,amscd,latexsym,mathrsfs,theorem}
\usepackage[all]{xy}

\newcommand{\F}{{\mathbb{F}}}

\newcommand{\N}{{\mathbb{N}}}

\newcommand{\adj}{\mathrm{adj}}

\newcommand{\id}{\mathrm{id}}
\newcommand{\into}{\hookrightarrow}

\newcommand{\Mod}{\mathbf{Mod}}
\newcommand{\Mor}{\mathrm{Mor}}

\newcommand{\spec}{\mathrm{Spec}\,}
\newcommand{\supp}{\mathrm{supp}\,}
\newcommand{\Aut}{\mathrm{Aut}}

\newcommand{\Gal}{\mathrm{Gal}}
\newcommand{\GL}{\mathrm{GL}\,}
\newcommand{\Hom}{\mathrm{Hom}}

\newcommand{\Iso}{\mathrm{Iso}\,}

\newcommand{\Ob}{\mathrm{Ob}\,}

\newcommand{\Rep}{\mathrm{Rep}\,}

\newcommand{\Sh}{{\mathcal S}}

\newcommand{\eo}{\mathfrak{o}}

\newcommand{\eE}{\mathcal{E}}

\newcommand{\eX}{{\mathfrak X}}

\newcommand{\ohne}{\setminus}

\newenvironment{theoremon}{\noindent {\bf Theorem \ref{main}}\it}{}
\newenvironment{proof}{\noindent {\bf Proof}}{\mbox{}\hspace*{\fill}$\Box$}
\theoremstyle{plain}

{\theorembodyfont{\slshape}

        \newtheorem{thm}{Theorem}[section]
        \newtheorem{cor}[thm]{Corollary}
        \newtheorem{lem}[thm]{Lemma}
        \newtheorem{prop}[thm]{Proposition}

}

{\theorembodyfont{\rmfamily}

        \newtheorem{defn}[thm]{Definition}
        
        \newtheorem{rem}[thm]{Remark}

}

\newcommand{\iproof}{{\bf Proof\ }}

\begin{document}
\title{Principal bundles on $p$-adic curves and parallel transport}
\author{Urs Hackstein}
\date{}
\maketitle

\begin{abstract}
We define functorial isomorphisms of parallel transport along \'etale paths for a class of principal $G$-bundles on a $p$-adic curve. Here $G$ is a connected reductive algebraic group of finite presentation and the considered principal bundles are just those with potentially strongly semistable reduction of degree zero. The constructed isomorphisms yield continous functors from the \'etale fundamental groupoid of the given curve to the category of topological spaces with a simply transitive continous right $G(\mathbb{C}_{p})$-action. This generalizes a construction in the case of vector bundles on a $p$-adic curve by Deninger and Werner. It may be viewed as a partial $p$-adic analogue of the classical theory by Ramanathan of principal bundles on compact Riemann surfaces, which generalizes the classical Narasimhan--Seshadri theory of vector bundles on compact Riemann surfaces.

\centerline{MSC (2000): 14H60, 14H30, 11G20}
\end{abstract}

\section*{Introduction}

Every finite dimensional complex representation of the fundamental group of a compact Riemann surface gives rise to a flat vector bundle and hence to a holomorphic vector bundle. By a classical result of Andr\'e Weil (\cite{W}), a holomorphic vector bundle on a compact Riemann surface is given by a representation of the fundamental group if and only if all its indecomposable components have degree zero. A famous theorem of Narasimhan and Seshadri states that there is an equivalence between the category of unitary representations of the fundamental group and those of polystable vector bundles of degree zero (\cite[Corollary 12.1]{Na-Se}). In particular, this result implies that every stable vector bundle of degree zero comes from an irreducible unitary representation (\cite[Corollary 12.2]{Na-Se}). These results have been extended to the case of principal $G$-bundles on a compact Riemann surface $X$ by Ramanathan in \cite{R}, where $G$ is a connected reductive algebraic group over $\mathbb{C}$. Ramanathan associates principal $G$-bundles to certain representations of \mbox{$\pi_{1}(X-x_{0})$} where $x_{0}$ is some point of $X$. Using an extension of Mumford's notion of semistability, he shows that a $G$-bundle is stable if and only if it is given by an irreducible unitary representation of $\pi_{1}(X-x_{0})$.

These results do not extend easily to the $p$-adic world as, for example, there is no notion of unitarity or universal coverings. However, Deninger and Werner constructed in their recent article \cite{De-We1} a partial $p$-adic analogue of the results by Narasimhan and Seshadri. They consider vector bundles $E$ on $X_{\mathbb{C}_{p}}$, where $X$ is a smooth projective curve over $\overline{\mathbb{Q}}_{p}$, which have potentially strongly semistable reduction of degree zero. For all such vector bundles $E$ they construct functorial isomorphisms of ``parallel transport'' along \'etale paths between the fibres of $E$. 
One obtains a representation $\rho_{E,x}$ of $\pi_{1}(X,x)$ on $x^{*}E=:E_{x}$ for every point $x\in X(\mathbb{C}_{p})$. Several aspects in relation to this construction have been studied by Deninger, Werner, Herz and Tong in \cite{De-We2}, \cite{He}, \cite{T} and \cite{We}. 

A natural question is whether one can extend the results of Deninger and Werner to the case of principal bundles in order to obtain a partial $p$-adic analogue of the classical results of Ramanathan. The construction of isomorphisms by Deninger and Werner may be extended easily, but the characterisation of the class of bundles to those they asoociate representations do not transfer directly to the case of principal bundles as they use frequently special properties of vector bundles, e.g. they work always in Zariski toplogy (whereas in the case of principal bundles fppf or \'etale topology is needed) or use the fact that vector bundles are coherent sheaves.

The present article, which is a shortened english version of the authors thesis (\cite{H}), shows now how the construction by Deninger and Werner can be extended to the case of principal $G$-bundles on $X_{\mathbb{C}_{p}}$, where $G$ is a connected reductive group scheme of finite presentation over the ring of integers $\mathfrak{o}$ of $\mathbb{C}_{p}$. Here and in the whole paper, a principal $G$-bundle or $G$-bundle for short is a $G$-torsor with respect to the fppf-topology and we consider $G$-bundles $E$ with potentially strongly semistable reduction of degree zero. We use here the notion of a degree of a principal bundle of Holla and Narasimhan (\cite{HN}) and the definition of semistability by Ramanathan (\cite{R}). As in the case of vector bundles we say that a $G$-bundle $E$ on $X_{\mathbb{C}_{p}}$ has \emph{strongly semistable reduction of degree zero} if $E$ extends to a $G$-bundle $\widetilde{E}$ on $\mathfrak{X}_{\mathfrak{o}}=\mathfrak{X}\otimes\mathfrak{o}$, for a suitable model $\mathfrak{X}$ of $X$ such that the pullback of the special fibre $\widetilde{E}_{k}$ of $\widetilde{E}$ to the normalisation of every irreducible component of $\mathfrak{X}_{k}$ is strongly semistable of degree zero. Furthermore, we say that $E$ has \emph{potentially strongly semistable reduction of degree zero} if there is a finite \'etale morphism $\alpha\colon Y\to X$ such that $\alpha_{\mathbb{C}_{p}}^{*}E$ has strongly semistable reduction of degree zero. If we denote by $\mathcal{P}(G(\mathbb{C}_{p}))$ the category of topological spaces with a simply transitive continous $G(\mathbb{C}_{p})$-action from the right, our main result is the following:

\begin{theoremon}
 Let $E$ be a pincipal $G$-bundle on $X_{\mathbb{C}_{p}}$ with potentially strongly semi\-stable reduction of degree zero. Then there exist functorial isomorphisms of ``parallel transport'' along \'etale paths of the fibres of $E$ on $X_{\mathbb{C}_{p}}$. In particular, for every such principal $G$-bundle $E$, one gets a continuous functor $\rho_{E}$ from $\pi_{1}(X)$ to $\mathcal{P}(G(\mathbb{C}_{p}))$ such that $\rho_{E}(x)=E_{x}=x^{*}E$ for all $x\in X(\mathbb{C}_{p})$. 
This defines a functor \[\rho\colon\mathscr{B}_{X_{\mathbb{C}_{p}}}^{ps}(G)\to\Rep_{\pi_{1}(X)}(G(\mathbb{C}_{p})),\, E\mapsto\rho_{E},\] where $\mathscr{B}_{X_{\mathbb{C}_{p}}}^{ps}(G)$ is the category of $G$-bundles with potentially strongly semistable reduction of degree zero and $\Rep_{\pi_{1}(X)}(G(\mathbb{C}_{p}))$ the category of continous functors from $\pi_{1}(X)$ to $\mathcal{P}(G(\mathbb{C}_{p}))$.
\end{theoremon}

This construction may be considered like a partial $p$-adic analogue of the theory of Ramanathan.

The paper is organized as follows. 

In the first section we recall some basic facts about torsors and the \'etale fundamental groupoid. Furthermore, we recall the definitions of certain categories of coverings defined in \cite{De-We1} which we need for the rest of the paper.

In the second section we define and investigate categories $\mathscr{B}_{X_{\mathbb{C}_{p}},D}(G)$ and $\mathscr{B}_{\mathfrak{X}_{\mathfrak{o}},D}(G)$ involving a divisor $D$ on $X$ where $\mathfrak{X}$ is a fixed model of $X$. 
For every bundle in these two categories, we construct isomorphisms of parallel transport along \'etale paths in $U = X \ohne D$ based on initial ideas of Deninger (\cite{De-We3}). First we do this for the bundles in the category $\mathscr{B}_{\mathfrak{X}_{\mathfrak{o}}, D}(G)$ for a fixed model $\mathfrak{X}$ of $X$ by imitating the construction of parallel transport in the case of vector bundles. Then we extend this construction to all bundles in the category $\mathscr{B}_{X_{\mathbb{C}_{p}},D}(G)$. In more technical terms, we construct a continuous functor $\rho_{E}$ from $\pi_{1}(X-D)$ to $\mathcal{P}(G(\mathbb{C}_{p}))$ (resp. $\mathcal{P}(G(\mathfrak{o}))$) such that $\rho_{E}(x)=E_{x}=x^{*}E$ for all $x\in X(\mathbb{C}_{p})$, where $\mathcal{P}(G(\mathfrak{o}))$ is the category of topological spaces with a simply transitive continous $G(\mathfrak{o})$-action from the right. By using a Seifert--van Kampen theorem, we show that the functors $\rho_{E}$ glue for different choices of the divisors $D$ (Proposition \ref{seifert}). We conclude this section by a collection of nice functorial properties of the construction and show that it is compatible with the construction in the case of vector bundles (Propositions \ref{pblvbl}, \ref{vb1}).

In the third section, we investigate the categories $\mathscr{B}_{\mathfrak{X}_{\mathfrak{o}},D}(G)$ and $\mathscr{B}_{X_{\mathbb{C}_{p}},D}(G)$ more carefully. First we prove that for a $G$-bundle $E$ on $\mathfrak{X}_{\mathfrak{o}}$ to be contained in $\mathscr{B}_{\eX_{\eo , D}}(G)$ it suffices that $\pi^*_k \eE_k$ is trivial, where $\pi_k$ is the special fibre of some $\pi$ as above. We continue with a characterisation of the principal bundles on a purely one-dimensional proper scheme over a finite field $\F_q$ whose pullback to the normalisation of each irreducible component is strongly semistable of degree zero: these are exactly the bundles whose pullback by a finite surjective morphism to a purely one-dimensional proper $\F_q$-scheme becomes trivial. For principal bundles on smooth projective curves over finite fields this characterisation follows from a result of Deligne (\cite{Las}). 

In the fourth and last section we finally state and prove our main theorems that our previous constructions together with the characterisation of the categories in the third section yield us functorial isomorphisms of parallel transport for all principal $G$-bundles with (potentially) strongly semistable reduction of degree zero and show that the fibr functor in a fixed point is faithful. We conclude with a collection of remaining open questions.

{\bf Acknowledgements:} I am indebted to Christopher Deninger for his continuous excellent support and encouragement during my work on my Phd-thesis. I am also grateful to Yves Laszlo for drawing my attention to a result of Deligne stated in \cite{Las}. Moreover, I would like to thank among many others Pham Ngoc Duy, Philippe Gille, Sylvain Maugeais, Michel Raynaud, Jean-Marc Fontaine, Alexander Steinmetz, Gabriel Herz, Stefan Wiech, Thomas Ludsteck and Jilong Tong for interesting and very helpful remarks and discussions. Some crucial parts of the present article have been investigated during a stay at the Laboratoire de Math\'ematiques of the Universit\'e de Paris-Sud at Orsay which was financed by the European Network Arithmetic Algebraic Geometry (MRTN-CT-2003-504917). I thank the work group ``Arithm\'etique et g\'eom\'etrie alg\'ebrique'' for its hospitality and in particular the local coordinator of the network, \'Etienne Fouvry, for excellent administrative support during my stay. Furthermore I was also partially supported by the Deutsche Forschungsgemeinschaft via its SFB 478 ``Geometrische Strukturen in der Mathematik''.

\section{Preliminaries}

\subsection{Torsors}
We start by recalling some basic definitions and facts about torsors.

Let $G$ be a group scheme which is faithfully flat and locally of finite presentation over some base scheme $S$. We call an $S$-scheme $P$ on which $G$ acts from the right via $P\times G\to P, (x,g)\mapsto xg$ a \emph{(right) $G$-torsor} with respect to the fppf topology if there is a covering $(U_{i}\to S)_{i\in I}$ in the fppf topology on $S$ such that $P\times_{S}U_{i}$ endowed with the induced $G\times_{S}U_{i}$-action is isomorphic to $G\times_{S}U_{i}$, endowed with the $G\times_{S}U_{i}$-action induced by the group multiplication. In the following, a \emph{$G$-torsor} is always a (right) $G$-torsor with respect to the fppf topology. A $G$-torsor $P$ is called \emph{trivial} if $P$ is isomorphic to $G$, considered as a $G$-torsor with right multiplication. It is well known that an $S$-scheme $P$ with a $G$-action from the right is a $G$-torsor if and only if $P$ is faithfully flat and locally of finite presentation over $S$ and the canonical morphism \mbox{$P\times G\to P\times P,$} \mbox{$(x,g)\mapsto (x,xg)$} is an isomorphism.

If $G$ is  a group object in the category of sheaves with respect to the fppf topology on $S$, then a right-$G$-(sheaf) torsor $P$ is defined as a sheaf $P$ on which $G$ acts by a morphism $P\times_{S}G\to P, (x,g)\mapsto xg$ such that the following condition holds.
There exists a covering $(U_{i}\to S)_{i\in I}$ for the fppf topology on $S$ such that $P\times_{S}U_{i}$ is isomorphic to $G\times_{S}U_{i}$ and this isomorphism is compatible with the actions of $G\times_{S}U_{i}$.
We call a $G$-torsor $P$ \emph{trivial} if it is isomorphic to $G$ as a $G$-torsor. It is obvious that a sheaf torsor is a torsor if it is representable by a scheme. It follows from Yoneda's Lemma that two torsors $P$ and $P'$ are isomorphic if and only if they are as sheaf torsors. Not every sheaf torsor is representable by a scheme, but in the case that $G$ is affine - like it will be in the present paper - all sheaf torsors are representable.

More generally, torsors can be defined similary for arbitrary Grothendieck topologies and arbitrary topoi. We conclude this section by recalling that there is a well known connection between torsors and cohomology:

\begin{prop}
If $G$ is a flat group scheme of finite type over a scheme $X$, then we can identify the set of isomorphism classes of $G$-sheaf-torsors with respect to the fppf-topology (resp. Zariski-topology resp. \'etale topology) on $X$ with $\check{H}^{1}(X_{fppf},G)$ (resp. $\check{H}^{1}(X,G)$ resp. $\check{H}^{1}_{\acute{e}t}(S,G)$).\\
If $G$ is in addition affine, then we can identify the set of isomorphism classes of $G$-torsors with respect to the fppf topology (resp. Zariski topology resp. \'etale topology) on $X$ with $\check{H}^{1}(X_{fppf},G)$ (resp. $ \check{H}^{1}(X,G)$ resp. $\check{H}^{1}_{\acute{e}t}(S,G)$).
\end{prop}

\subsection{Categories of coverings}
In this section, we recall the definition and some basic properties of some categories of coverings introduced by Deninger and Werner in \cite{De-We1}. In the following,
a variety over a field $k$ is always a geometrically irreducible and geometrically reduced and separated scheme of finite type over $k$. A curve is a one-dimensional variety.\par
We let $R$ be a valuation ring with quotient field $Q$ of characteristic zero. For a smooth projective curve $X$ over $Q$, we consider a model $\mathfrak{X}$ of $X$ over $R$, i.\,e. a finitely presented, flat and proper scheme $\mathfrak{X}$ over $\spec\,R$ such that $X$ is isomorphic to $\mathfrak{X}\otimes_{R}Q$. Finally, we write $X-D$ for $X-\supp\,D$ for a divisor $D$ of $X$. In this situation we may define:

\begin{defn}[\protect{\cite[p. 4,5]{De-We1}}]
We define a category $\Sh_{\mathfrak{X},D}$ as follows:
The objects are finitely presented, proper $R$-morphism \mbox{$\pi\colon\mathcal{Y}\to\mathfrak{X}$} whose generic fibre \mbox{$\pi_{Q}\colon\mathcal{Y}_{Q}\to X$} is finite and such that the induced morphism $\pi_{Q}\colon\pi_{Q}^{-1}(X-D)\to X-D$ is \'etale.
A morphism from $\pi_{1}\colon\mathcal{Y}_{1}\to\mathfrak{X}$ to $\pi_{2}\colon\mathcal{Y}_{2}\to\mathfrak{X}$ is a morphism $\varphi\colon\mathcal{Y}_{1}\to\mathcal{Y}_{2}$ such that $\pi_{1}=\pi_{2}\circ\varphi$.
If such a morphism exists, we say that $\pi_{1}$ \emph{dominates} $\pi_{2}$. If, in addition, $\varphi_{Q}=\varphi\otimes_{R}Q$ induces an isomorphism of the local rings of two generic points, we say that $\pi_{1}$ \emph{strictly dominates} $\pi_{2}$.
\end{defn}

It is obvious that finite products and finite fibre products exist in $\Sh_{\mathfrak{X},D}$. Furthermore, it is easy to verify that for every morphism $\varphi\colon(\pi_{1}\colon\mathcal{Y}_{1}\to\mathfrak{X})\to(\pi_{2}\colon\mathcal{Y}_{2}\to\mathfrak{X})$ the induced morphism $\varphi\colon\mathcal{Y}_{1}\to\mathcal{Y}_{2}$ is proper and of finite type and its generic fibre $\varphi_{Q}$ is finite and \'etale over $X-D$.

We recall also the full subcategory $S_{\mathfrak{X},D}^{\text{good}}$ of $\Sh_{\mathfrak{X},D}$:

\begin{defn}
The full subcategory $\Sh_{\mathfrak{X},D}^{\text{good}}$ of $\Sh_{\mathfrak{X},D}$ consists of all elements \mbox{$\pi\colon\mathcal{Y}\to\mathfrak{X}$} whose structure morphism $\lambda\colon\mathcal{Y}\to\spec\,R$ is flat and satisfies \mbox{$\lambda_{*}\mathcal{O}_{\mathcal{Y}}=\mathcal{O}_{\spec\,R}$} universally, and whose generic fibre \mbox{$\lambda_{Q}\colon\mathcal{Y}_{Q}\to\spec\,Q$} is smooth.
\end{defn}

It is easy to see that if $\pi\colon\mathcal{Y}\to\mathfrak{X}$ lies in $\Sh_{\mathfrak{X},D}^{\text{good}}$, then $\mathcal{Y}_{Q}$ is geometrically connected and hence a smooth projective curve over $Q$. This implies that $\mathcal{Y}$ is irreducible and reduced.

As a corollary of an analogous statement for abitrary discrete valuation rings, Denin\-ger and Werner proved the following result: 

\begin{prop}\label{cor1good}
Let $X$ be a smooth projective curve over $\overline{\mathbb{Q}}_{p}$, $D$ a divisor in $X$ and $\mathfrak{X}$ a model of $X$ over $\spec\,\overline{\mathbb{Z}}_{p}$. Furthermore, let $\pi_{i}\colon\mathcal{Y}_{i}\to\mathfrak{X}$ be a finite number of objects in $\Sh_{\mathfrak{X},D}$.
Then there exists a finite extension $K$ of $\mathbb{Q}_{p}$ and a smooth projective curve $X_{K}$ over $K$ with a model $\mathfrak{X}_{\mathfrak{o}_{K}}$ over $\mathfrak{o}_{K}$ and a divisor $D_{K}$ of $X_{K}$ such that the following holds: We have $X=X_{K}\otimes_{K}\overline{K}$, $D=D_{K}\otimes_{K}\overline{K}$ and $\mathfrak{X}=\mathfrak{X}_{\mathfrak{o}_{K}}\otimes_{\mathfrak{o}_{K}}\overline{\mathbb{Z}}_{p}$. Furthermore there is a element $\pi_{\mathfrak{o}_{K}}\colon\mathcal{Y}_{\mathfrak{o}_{K}}\to\mathfrak{X}_{\mathfrak{o}_{K}}$ of $\Sh_{\mathfrak{X}_{\mathfrak{o}_{K}},D_{K}}^{\text{good}}$ such that the morphism $\pi=\pi_{\mathfrak{o}_{K}}\otimes_{\mathfrak{o}_{K}}\overline{\mathbb{Z}}_{p}$ dominates all $\pi_{i}$.
\end{prop}

\begin{cor}\label{cor2good}
Let $X, \mathfrak{X}$ and $D$ be as before. Then any finite number of objects \mbox{$\pi_{i}\colon\mathcal{Y}_{i}\to\mathfrak{X}$} in $S_{\mathfrak{X},D}$ is dominated by a common object $(\pi\colon\mathcal{Y}\to\mathfrak{X})\in S_{\mathfrak{X},D}^{\text{good}}$.
\end{cor}

Next we need another full subcategory of $\Sh_{\mathfrak{X},D}$:

\begin{defn}[\protect{\cite[p. 5]{De-We1}}]
The full subcategory $\Sh_{\mathfrak{X},D}^{\text{ss}}$ in $\Sh_{\mathfrak{X},D}$ consists of all elements $\pi\colon\mathcal{Y}\to\mathfrak{X}$ in $\Sh_{\mathfrak{X},D}$ such that $\lambda\colon\mathcal{Y}\to\spec\,R$ is a semistable curve whose generic fibre $\mathcal{Y}_{Q}$ is a smooth projective curve over $Q$. 
\end{defn}

\begin{rem}\label{doms}
We recall from \cite[Theorem 1, Corollary 3]{De-We1} that all stated facts on the subcategory $\Sh_{\mathfrak{X},D}^{\text{good}}$ hold verbatim for the subcategory $\Sh_{\mathfrak{X},D}^{\text{ss}}$.
\end{rem}

\subsection{\'Etale fundamental groupoid and fibre functors}

Last, but not least, we recall the notions of an \'etale fundamental groupoid and of fibre functors which will be important in the following. The general reference is \cite{SGA1}, but a collection of basic facts is also given in \cite[p.577 etc.]{De-We1}.

Let $Z$ be a variety over $\overline{\mathbb{Q}}_{p}$, and let $z\in Z(\mathbb{C}_{p})$ be a geometric point. The fibre functor $F_{z}$ from the category of finite \'etale coverings $Z'$ of $Z$ to finite sets associates to every such $Z'$ the set of its $\mathbb{C}_{p}$-valued points lying over $z$.

Moreover, we define the topological category $\pi_{1}(Z)$ as follows:
the objects are just the geometric points of $Z$, i.\,e. the set of objects is $Z(\mathbb{C}_{p})$.
A morphism in $\pi_{1}(Z)$ between two $\mathbb{C}_{p}$-valued points $z$ and $z^{*}$ is an isomorphism of its associated fibre functors. Such an isomorphism of fibre functors is called an \emph{\'etale path} from $z$ to $z^{*}$.

We know by \cite[Exposé V, 4., 7.]{SGA1} that every fibre functor is pro-representable. Therefore, $\Mor_{\pi_{1}(Z)}(z,z^{*})$ is a profinite set and a compact totally disconnected Hausdorff space (see \cite[p. 577]{De-We1}). Moreover, composition of morphisms induces a continous map \mbox{$\Mor_{\pi_{1}(Z)}(z,z^{*})\times\Mor_{\pi_{1}(Z)}(z^{*},z^{**})\to \Mor_{\pi_{1}(Z)}(z,z^{**})$}.
The topological category $\pi_{1}(Z)$ is called the \emph{\'etale fundamental groupoid of $Z$}.

\section{Parallel transport for principal bundles}
\subsection{Parallel transport for principal bundles on $\mathfrak{X}\otimes_{\overline{\mathbb{Z}}_{p}}\mathfrak{o}$ for a fixed model $\mathfrak{X}$ of $X$}

Let $\mathfrak{o}$ be the ring of integers in $\mathbb{C}_{p}$ and let $S:=\spec\,\mathfrak{o}$. Furthermore, let $G:=\spec\,A$  be a smooth affine group scheme of finite presentation over $\mathfrak{o}$. It follows directly that $G$ is also faithfully flat over $\mathfrak{o}$, since any flat group scheme of finite presentation is faithfully flat (\cite[Expos\'e $VI_{B}$, Proposition 9.2. (xi),(xii)]{SGA3}).

In the remaining part of this article, a $G$-torsor or \emph{principal $G$-bundle} $P$ on a $\mathfrak{o}$-scheme $\xi$ is always a representable right $G_{\xi}=G\times_{\spec\mathfrak{o}}\xi$-torsor with respect to the fppf topology on $\xi$. By faithfully flat descent, it follows that such a $G$-torsor on an $\mathfrak{o}$-scheme $\xi$ is always smooth over $\xi$. Furthermore, we deduce from \cite[Th\'eor\`eme 11.7.]{Gr} that $P$ is a $G$-torsor with respect to the \'etale topology.

Let now $P$ be a principal $G$-bundle on an $\mathfrak{o}$-scheme $\xi$. Then we have a natural action of $G(\mathfrak{o})$ on $\Gamma(\xi, P)$ which makes $\Gamma(\xi, P)$ into a right-$G(\mathfrak{o})$-module.

Let $f\colon\xi_{1}\to\xi_{2}$ be a morphism of $\mathfrak{o}$-schemes, and let $P$ be a $G$-bundle on $\xi_{2}$. Then the pullback $f^{*}P:=\xi_{1}\times_{\xi_{2}}P$ is a principal $G$-bundle on $\xi_{1}$, and there is a natural map \mbox{$f^{*}\colon\Gamma(\xi_{2},P)\to\Gamma(\xi_{1},f^{*}P)$} which is equivariant under the right-$G(\mathfrak{o})$ action.

For $n\geq 1$ we set $\mathfrak{o}_{n}:=\mathfrak{o}/p^{n}\mathfrak{o}=\overline{\mathbb{Z}}_{p}/p^{n}\overline{\mathbb{Z}}_{p}$ such that we have $\mathfrak{o}=\varprojlim\limits_{n}\mathfrak{o}_{n}$. Moreover, we define $A_{n}:=A\otimes_{\mathfrak{o}}\mathfrak{o}_{n}$ and denote by $G_{n}:=\spec\,A_{n}$ the reduction of $G$ modulo $p^{n}$. For any $\mathfrak{o}$-scheme $\xi$ we set $\xi_{n}:=\xi\otimes_{\mathfrak{o}}\mathfrak{o}_{n}$. For a $G$-bundle $P$ on $\xi$, we define the scheme $P_{n}$ as the reduction of $P$ to a $G$- resp. $G_{n}$-torsor on $\xi_{n}$, i.\,e. we set $P_{n}=i_{n}^{*}P$, where $i_{n}\colon\xi_{n}\to\xi$ is the canonical morphism.

We now endow $G(\mathfrak{o})$ with a prodiscrete topology induced by the canonical isomorphism \[G(\mathfrak{o})\cong\Hom_{\mathfrak{o}\text{-Alg.}}(A,\mathfrak{o})\cong\varprojlim_{n}\Hom_{\mathfrak{o}\text{-Alg.}}(A,\mathfrak{o}_{n})\cong\varprojlim_{n}G(\mathfrak{o}_{n})\] and the discrete topology on $G(\mathfrak{o}_{n})$ such that $G(\mathfrak{o})$ becomes a topological group. The formal smoothness of $G$ over $\mathfrak{o}$ implies 
that all transition maps $G(\mathfrak{o}_{n+1})\to G(\mathfrak{o}_{n})$ are surjective. From the construction of the  projective limit follows thus that all morphisms \mbox{$G(\mathfrak{o})\to G(\mathfrak{o}_{n})$} are also surjective.

These prerequisites allow us to define a category $\mathscr{B}_{\mathfrak{X}_{\mathfrak{o}},D}(G)$ of $G$-torsors on $\mathfrak{X}\otimes_{\overline{\mathbb{Z}}_{p}}\mathfrak{o}$ where $\mathfrak{X}$ is a fixed model of $X$ over $\overline{\mathbb{Z}}_{p}$. 

\begin{defn}\label{modcat}
We fix a model $\mathfrak{X}$ over $\overline{\mathbb{Z}}_{p}$ of a smooth projective curve $X$ over $\overline{\mathbb{Q}}_{p}$ and a divisor $D$ on $X$. We define $\mathscr{B}_{\mathfrak{X}_{\mathfrak{o}},D}(G)$ to be the category of principal $G$-bundles $P$ on $\mathfrak{X}_{\mathfrak{o}}=\mathfrak{X}\otimes_{\overline{\mathbb{Z}}_{p}}\mathfrak{o}$ with the following property:

For all $n\geq 1$ there exists an object $\pi\colon\mathcal{Y}\to\mathfrak{X}$ of $S_{\mathfrak{X},D}$ such that $\pi_{n}^{*}P_{n}$ is trivial, i.\,e. it is isomorphic as $G_{n}$-torsors to the trivial torsor $\mathcal{Y}_{n}\otimes_{\mathfrak{o}_{n}}G_{n}$ where $\pi_{n}\colon\mathcal{Y}_{n}\to\mathfrak{X}_{n}$ is the morphism induced by $\pi$.
\end{defn}

In the following we fix a smooth projective curve $X$ over $\overline{\mathbb{Q}}_{p}$, a model $\mathfrak{X}$ of $X$ over $\overline{\mathbb{Z}}_{p}$ and a divisor $D$ in $X$. We easily see that for a morphism \mbox{$f\colon\mathfrak{X}\to\mathfrak{X'}$} of models over $\overline{\mathbb{Z}}_{p}$ of smooth projective curves $X$ and $X'$ over $\overline{\mathbb{Q}}_{p}$ and a divisor $D'$ in $X'$, the pullback functor $f^{*}$ for $G$-torsors induces a functor \mbox{$f^{*}\colon\mathscr{B}_{\mathfrak{X'}_{\mathfrak{o}},D'}(G)\to\mathscr{B}_{\mathfrak{X}_{\mathfrak{o}},f^{*}D'}(G)$}.
 
Let $x\colon\spec\,\mathbb{C}_{p}\to X$ be a $\mathbb{C}_{p}$-valued point of $X$, and let $P$ be a principal $G$-bundle over $\mathfrak{X}_{\mathfrak{o}}$. Then $x$ induces a morphism $\widetilde{x}\colon\spec\,\mathbb{C}_{p}\to X\to\mathfrak{X}$ which again induces, by the properness of $\mathfrak{X}$, 
         an unique morphism $x_{\mathfrak{o}}\colon\spec\mathfrak{o}\to\mathfrak{X}_{\mathfrak{o}}$. For $n\in\N ,n\geq1$ we set \mbox{$x_{n}:=x_{\mathfrak{o}}\circ(\spec\,\mathfrak{o}_{n}\to\spec\,\mathfrak{o})$} such that we can define the groups $P_{x_{\mathfrak{o}}}$ and $P_{x_{n}}$ by $P_{x_{\mathfrak{o}}}:=(x_{\mathfrak{o}}^{*}P)(\mathfrak{o})$ and \mbox{$P_{x_{n}}:=(x_{n}^{*}P)(\mathfrak{o}_{n})$}. Next we want to define topologies on $P_{x_{\mathfrak{o}}}$ and $P_{x_{n}}$. We use here that the $G$-torsor $x_{\mathfrak{o}}^{*}P$ is affine by descent theory. As $(x_{\mathfrak{o}}^{*}P)_{n}=x_{n}^{*}P$, we have the following lemma:

\begin{lem}
We have $P_{x_{\mathfrak{o}}}=\varprojlim\limits_{n}P_{x_{n}}$.
\end{lem}
\begin{proof}{ As $x_{\mathfrak{o}}^{*}P$ is affine over $\mathfrak{o}$, say \mbox{$x_{\mathfrak{o}}^{*}P=\spec\,B$}, we have that $P_{x_{\mathfrak{o}}}=(x_{\mathfrak{o}}^{*}P)(\mathfrak{o})$ is isomorphic to \mbox{$\Hom_{\mathfrak{o}\text{-Alg}}(B,\mathfrak{o})=\varprojlim\limits_{n}\Hom_{\mathfrak{o}\text{-Alg}}(B,\mathfrak{o}_{n})$.} Moreover, $f(p^{n}B)=0$ holds for every morphism of $\mathfrak{o}$-algebras $f\colon B\to\mathfrak{o}_{n}$. Every such morphism $f$ factorizes therefore through $B/p^{n}B=B\otimes_{\mathfrak{o}}\mathfrak{o}_{n}$. This implies \mbox{$\Hom_{\mathfrak{o}\text{-Alg}}(B,\mathfrak{o}_{n})\cong \Hom_{\mathfrak{o}_{n}\text{-Alg}}(B\otimes_{\mathfrak{o}}\mathfrak{o}_{n},\mathfrak{o}_{n})$}. Using $(x_{\mathfrak{o}}^{*}P)_{n}=x_{n}^{*}P$ we obtain thus:
\begin{align*}
P_{x_{\mathfrak{o}}} & \cong\varprojlim_{n}Hom_{\mathfrak{0}_{n}-Alg}(B\otimes_{\mathfrak{o}}\mathfrak{o}_{n},\mathfrak{o}_{n})\\
& \cong\varprojlim_{n}(i_{n}^{*}(x_{\mathfrak{o}}^{*}P))(\mathfrak{o}_{n})\cong\varprojlim_{n}(x_{n}^{*}P)(\mathfrak{o}_{n})=\varprojlim_{n}P_{x_{n}}
\end{align*}}\end{proof}

If we endow $P_{x_{n}}$ with the discrete topology for all $n$, then we obtain a prodiscrete topology on $P_{x_{\mathfrak{o}}}$. This makes $P_{x_{\mathfrak{o}}}$ into a topological space with a continous $G(\mathfrak{o})$-action. 
Furthermore, we deduce from the fact that $\mathfrak{o}$ is a complete henselian ring that $x_{\mathfrak{o}}^{*}P$ is a trivial $G$-torsor with respect to the \'etale (resp. fppf) topology. Thus $G(\mathfrak{o})$ and $G(\mathfrak{o}_{n})$ act simply transitive on $P_{x_{\mathfrak{o}}}$ resp. $P_{x_{n}}$.

Next we deduce easily from descent theory that $x_{\mathfrak{o}}^{*}P$ is smooth over $\spec\,\mathfrak{o}$. This implies 
that all transition maps \mbox{$x_{\mathfrak{o}}^{*}P(\mathfrak{o}_{n+1})\to x_{\mathfrak{o}}^{*}P (\mathfrak{o_{n}})$} and, by the construction of the projective limit, all morphisms $x_{\mathfrak{o}}^{*}P(\mathfrak{o})\to x_{\mathfrak{o}}^{*}P(\mathfrak{o}_{n})$ are surjective.

\begin{defn}
We define the category $\mathcal{P}(G(\mathfrak{o}_{n}))$ for every $n\in\N, n\geq 1$ as follows:
The objects are sets with a simply transitive right $G(\mathfrak{o}_{n})$-action and its morphisms are just the $G(\mathfrak{o}_{n})$-equivariant maps between them.
\end{defn}

If we fix now two arbitrary topological spaces $P_{1},P_{2}\in\mathcal{P}(G(\mathfrak{o}))$, then the choice of two elements $p_{1}\in P_{1}$ and $p_{2}\in P_{2}$ defines a bijection \[\Phi\colon G(\mathfrak{o})\xrightarrow{\sim}\Mor_{\mathcal{P}(G(\mathfrak{o}))}(P_{1},P_{2}),\,
g \longmapsto(\varphi_{g}\colon P_{1}\to P_{2},\,p_{1}h\mapsto p_{2}gh)\] which follows from the fact that $G(\mathfrak{o})$ acts simply transitive continous on $P_{1}$ and $P_{2}$. One checks easily that the topology on $\Mor_{\mathcal{P}(G(\mathfrak{o}))}(P_{1},P_{2})$ making it into a homeomorphism does not depend on the choice of $p_{1}$ and $p_{2}$. Thus $\mathcal{P}(G(\mathfrak{o}))$ becomes in this way a topological category.

Let now $P$ be a $G$-torsor in $\mathscr{B}_{\mathfrak{X}_{\mathfrak{o}},D}(G)$. Then we choose for all $n\geq 1$ an object $(\pi\colon\mathcal{Y}\to\mathfrak{X})\in\Sh_{\mathfrak{X},D}^{\text{good}}$ such that $\pi_{n}^{*}P_{n}$ is trivial on $\mathcal{Y}_{n}$.
We set $U=X-D$ and $V=\mathcal{Y}\otimes_{\overline{\mathbb{Z}}_{p}}\overline{\mathbb{Q}}_{p}-\pi^{*}D$. Then $V$ is a finite \'etale covering of $U$. Thus there exists a point $y\in V(\mathbb{C}_{p})$ over $x\in U(\mathbb{C}_{p})$ by \cite[4.2]{Mu} and it holds:

\begin{lem}\label{20}
The pullback map $y_{n}^{*}\colon\Gamma(\mathcal{Y}_{n},\pi_{n}^{*}P_{n})\to\Gamma(\spec\,\mathfrak{o}_{n},y_{n}^{*}\pi_{n}^{*}P_{n})=P_{x_{n}}$ is an isomorphism of right-$G(\mathfrak{o})$ sets. 
\end{lem}
\begin{proof}{ It suffices to show that the map \[y_{n}^{*}\colon\Gamma(\mathcal{Y}_{n},\mathcal{Y}_{n}\times_{\spec\,\mathfrak{o}_{n}}G_{n})\to\Gamma(\spec\,\mathfrak{o}_{n},y_{n}^{*}(\mathcal{Y}_{n}\times_{\spec\,\mathfrak{o}_{n}}G_{n}))=\Gamma(\spec\,\mathfrak{o}_{n},G_{n})=G_{n}(\mathfrak{o}_{n})\] is bijective.

By the universal property of fibre products we have the following natural bijections:
\[\Gamma(\mathcal{Y}_{n},\mathcal{Y}_{n}\times_{\spec\,\mathfrak{o}_{n}}G_{n})\cong \Mor_{\spec\,\mathfrak{o}_{n}}(\mathcal{Y}_{n},G_{n})\cong\Hom_{\mathfrak{o}_{n}\text{-Alg}}(A_{n},\Gamma(\mathcal{Y}_{n},\mathcal{O}_{\mathcal{Y}_{n}}))\]
As $\pi\in\Sh_{\mathfrak{X},D}^{\text{good}}$ and thus \mbox{$\Gamma(\mathcal{Y}_{n},\mathcal{O}_{\mathcal{Y}_{n}})=(\lambda_{n})_{*}\mathcal{O}_{\mathcal{Y}_{n}}(\spec\,\mathfrak{o}_{n})=\mathcal{O}_{\spec\,\mathfrak{o}_{n}}(\spec\,\mathfrak{o}_{n})=\mathfrak{o}_{n}$} holds, this is in canonical way isomorphic to $\Hom_{\mathfrak{o}_{n}-Alg}(A_{n},\mathfrak{o}_{n})=G_{n}(\mathfrak{o}_{n})$. If we identify $\Gamma(\mathcal{Y}_{n},\mathcal{Y}_{n}\times_{\spec\,\mathfrak{o}_{n}}G_{n})$ with $G_{n}(\mathfrak{o}_{n})$, then the map $y_{n}^{*}$ is the identity.}\end{proof}\par

Lemma \ref{20} enables us to construct a functor $\rho_{P,n}\colon\pi_{1}(U)\to\mathcal{P}(G(\mathfrak{o}))$ for every \mbox{$P\in\mathscr{B}_{\mathfrak{X}_{\mathfrak{o}},D}(G)$} and all $n\in\N, n\geq 1$ using the notation $U=X-D$ like above:
 
For every $\mathbb{C}_{p}$-valued point $x$ of $U$ we set $\rho_{P,n}(x):=P_{x_{n}}$. An \'etale path $\gamma$ from $x\in U(\mathbb{C}_{p})$ to $x'\in U(\mathbb{C}_{p})$ is mapped to a morphism $\rho_{P, n}(\gamma)$ in $\mathcal{P}(G(\mathfrak{o}))$ as follows:

First we choose an object $(\pi\colon\mathcal{Y}\to\mathfrak{X})\in\Sh_{\mathfrak{X},D}^{\text{good}}$ such that $\pi_{n}^{*}P_{n}$ is trivial on $\mathcal{Y}_{n}$. Furthermore we choose a point $y\in V(\mathbb{C}_{p})$ over $x\in U(\mathbb{C}_{p})$. Then $\gamma y$ is a $\mathbb{C}_{p}$-valued point of $V$ over $x'$ where we denote by $\gamma y$ the image of $y$ under the map $\gamma_{V}\colon F_{x}(V)\to F_{x'}(V)$. Whithin these choices we are able to define $\rho_{P,n}(\gamma)$ as the composition \[\rho_{P,n}(\gamma)\colon P_{x_{n}}\xrightarrow{(y_{n}^{*})^{-1}}\Gamma(\mathcal{Y}_{n},\pi_{n}^{*}G_{n})\xrightarrow{(\gamma y)_{n}^{*}} P_{x'_{n}}.\]

Formally nearly in the same way as in the case of vector bundles (see \cite[Section 3]{De-We1}) one shows moreover that the definition of the morphism $\rho_{P,n}(\gamma)$ is independent of all our choices and defines a morphism in $\mathcal{P}(G(\mathfrak{o}_{n}))$. We obtain thus a well-defined continous functor $\rho_{P,n}\colon\pi_{1}(X-D)\to\mathcal{P}(G(\mathfrak{o}_{n}))$.

Using these functors $\rho_{P,n}$ we are now able to define a functor $\rho_{P}\colon\pi_{1}(U)\to\mathcal{P}(G(\mathfrak{o}))$:
For every $\mathbb{C}_{p}$-valued point $x$ of $U$ we set $\rho_{P}(x)=P_{x_{\mathfrak{o}}}$ and  for an \'etale path $\gamma$ we define $\rho_{P}(\gamma)$ as $\rho_{P}(\gamma)=\varprojlim_{n}\rho_{P,n}(\gamma)\colon P_{x_{\mathfrak{o}}}\to P_{x'_{\mathfrak{o}}}$.

This construction gives us a well defined functor $\rho_{P}$ because one verifies again formally nearly in the same  way as in the case of vector bundles (see the proof of \cite[Theorem 22]{De-We1}) that the family of maps $(\rho_{P,n}(\gamma))_{n\geq 1}$ defines a morphism between the projective systems $(P_{x_{n}})_{n\geq 1}$ and $(P_{x'_{n}})_{n\geq 1}$.

We should remark that for any fixed point $x\in U(\mathbb{C}_{p})$ we obtain a continous homomorphism $\rho_{P}\colon\pi_{1}(U,x)\to\Aut_{\mathcal{P}(G(\mathfrak{o}))}(P_{x_{\mathfrak{o}}})$.

If we denote by $\Rep_{\pi_{1}(U)}(G(\mathfrak{o}))$ the category of continous functors from $\pi_{1}(U)$ to $\mathcal{P}(G(\mathfrak{o}))$, then we obtain in conclusion a well-defined functor \[\rho\colon\mathscr{B}_{\mathfrak{X}_{\mathfrak{o},D}}(G)\to \Rep_{\pi_{1}(U)}(G(\mathfrak{o}))\] by mapping an object $P$ of $\mathscr{B}_{\mathfrak{X}_{\mathfrak{o},D}}(G)$ onto $\rho_{P}$ and a morphism \mbox{$f\colon P\to P'$} onto the family $(f_{x_{\mathfrak{o}}})_{x\in U(\mathbb{C}_{p})}$. Here $f_{x_{\mathfrak{o}}}\colon P_{x_{\mathfrak{o}}}\to P'_{x_{\mathfrak{o}}}$ denotes the map $f_{x_{\mathfrak{o}}}=(x_{\mathfrak{o}}^{*}f)(\mathfrak{o})$.

This functor is well-defined because the family of morphisms of topological spaces with a continous $G(\mathfrak{o})$-action $(f_{x_{\mathfrak{o}}}=x_{\mathfrak{o}}^{*}f\colon P_{x_{\mathfrak{o}}}\to P'_{x_{\mathfrak{o}}})_{x\in U(\mathbb{C}_{p})=\Ob\pi_{1}(U)}$ defines a natural transformation $\rho_{f}$ from $\rho_{P}$ to $\rho_{P'}$ which may be checked easily.

\subsection{Extension to principal bundles over $X_{\mathbb{C}_{p}}$}

Of course, our real interest is with principal bundles over $X_{\mathbb{C}_{p}}$ and not on a particular model. Our aim is to define functorial isomorphisms for principal bundles over $X_{\mathbb{C}_{p}}$ based on the isomorphisms of parallel transport defined for bundles on $\mathfrak{X}_{\mathfrak{o}}$ for a particular model $\mathfrak{X}$ in the previous section. This motivates the following definition:

\begin{defn}
For a smooth affine group scheme $G$ of finite presentation over $\mathfrak{o}$ let $\mathscr{B}_{X_{\mathbb{C}_{p}}, D}(G)$ be the full subcategory of all right $G$-torsors (i.\,e. $G_{X_{\mathbb{C}_{p}}}$-torsors) $P$ on $X_{\mathbb{C}_{p}}$ with the following property:\\
There is a model $\mathfrak{X}$ of $X$ over $\overline{\mathbb{Z}}_{p}$ and a $G$-torsor $\widetilde{P}$ on $\mathfrak{X}_{\mathfrak{o}}$ which belongs to the category $\mathscr{B}_{\mathfrak{X}_{\mathfrak{o}},D}(G)$ such that the $G$-torsors $P$ and $j^{*}\widetilde{P}$ on $X_{\mathbb{C}_{p}}$ are isomorphic where \mbox{$j\colon X_{\mathbb{C}_{p}}\to\mathfrak{X}_{\mathfrak{o}}$} is the natural immersion.
\end{defn}

Using the same arguments as for vector bundles, one defines starting from the construction in Section §2.1 for a fixed model a functor \[\rho\colon\mathscr{B}_{X_{\mathbb{C}_{p}}, D}(G)\to \Rep_{\pi_{1}(U)}(G(\mathbb{C}_{p}))\] 
into the category of continous functors from $\pi_{1}(U)$ into the category $\mathcal{P}(G(\mathbb{C}_{p}))$ of topological spaces with a simply transitive continous right $G(\mathbb{C}_{p})$-action:

For any torsor $P$ in $\mathscr{B}_{X_{\mathbb{C}_{p}}, D}(G)$, we define a continous functor \[\rho(P)=\rho_{P}\colon\pi_{1}(U)\to\mathcal{P}(G(\mathbb{C}_{p}))\] as follows:
For $x\in U(\mathbb{C}_{p})$, we set $\rho_{P}(x)=P_{x}=x^{*}P$. For two points \mbox{$x,x'\in U(\mathbb{C}_{p})$}, we define a continous map  \mbox{$\rho_{P}=\rho_{P,x,x'}\colon\Mor_{\pi_{1}(U)}(x,x')\to\Hom_{\mathbb{C}_{p}}(P_{x},P_{x'})$} by setting $\rho_{P}(\gamma)=\psi^{-1}_{x'}\circ(\rho_{\widetilde{P}}(\gamma)\otimes_{\mathfrak{o}}\mathbb{C}_{p})\circ\psi_{x}$. Here $\mathfrak{X}$ is a choosen model of $X$ over $\overline{\mathbb{Z}}_{p}$ and $\widetilde{P}$ a $G$-torsor in $\mathscr{B}_{\mathfrak{X}_{\mathfrak{o}},D}(G)$ together with an isomorphism \mbox{$\psi\colon P\to j^{*}\widetilde{P}$} of $G$-torsors over $X_{\mathbb{C}_{p}}$ where $j\colon X_{\mathbb{C}_{p}}\to\mathfrak{X}_{\mathfrak{o}}$ is the natural immersion. Furthermore, $\psi_{x}$ is the fibre map $\psi_{x}=x^{*}(\psi)\colon P_{x}\xrightarrow{\sim}(j^{*}\widetilde{P})_{x}=\widetilde{P}_{x_{\mathfrak{o}}}\otimes_{\mathfrak{o}}\mathbb{C}_{p}=\widetilde{P}_{x_{\mathfrak{o}}}\otimes_{\mathbb{Z}}\mathbb{Q}$. On the other hand, a morphism $f\colon P\to P'$ in $\mathscr{B}_{\mathfrak{X}_{\mathbb{C}_{p}},D}(G)$ is mapped onto the morphism $\rho(f)=\rho_{f}\colon\rho_{P}\to\rho_{P'}$ given by the family of maps \mbox{$(f_{x}=x^{*}(f)\colon P_{x}\to P_{x'})_{x\in U(\mathbb{C}_{p})}$}.

\begin{prop}
The above definition of the functor $\rho_{P}$ is independent of all our choices, i.\,e. of the choice of a model $\mathfrak{X}$ and the $G$-torsor $\widetilde{P}\in\mathscr{B}_{\mathfrak{X}_{\mathfrak{o},D}}(G)$.
\end{prop}

This follows from analogous results to \cite[Propositions 26, 27]{De-We1}.

\subsection{Gluing of the functors $\rho_{P}$}

If $P$ is a $G$-torsor on $X_{\mathbb{C}_{p}}$ which belongs to categories $\mathscr{B}_{\mathfrak{X}_{\mathfrak{o}},D_{1}}(G)$ and $\mathscr{B}_{\mathfrak{X}_{\mathfrak{o}},D_{2}}(G)$
for different divisors $D_{1}$ and $D_{2}$, we deduce from the Seifert--van Kampen theorem \cite[Proposition 34]{De-We1} that we can glue the functors $\rho_{1,P}$ and $\rho_{2.P}$ constructed with respect to the divisor $D_{1}$ resp. $D_{2}$:

\begin{prop}\label{seifert}
Let $D_{1}$ and $D_{2}$ be two divisors in $X$ and set $U_{1}=X-D_{1}$ and \mbox{$U_{2}=X-D_{2}$}. Let $P$ be a $G$-torsor on $X$ which belongs as to $\mathscr{B}_{X_{\mathbb{C}_{p}}, D_{1}}(G)$ as to $\mathscr{B}_{X_{\mathbb{C}_{p}}, D_{2}}(G)$ and let $\rho_{P}^{i}$ be our constructed continous functors \mbox{$\pi_{1}(U_{i})\to\mathcal{P}(G(\mathbb{C}_{p}))$}, where \mbox{$i=1,2$}.\\
Then there exists an unique continous functor $\rho_{P}\colon\pi(U_{1}\cup U_{2})\to\mathcal{P}(G(\mathbb{C}_{p}))$ which induces the functors $\rho_{P}^{i}$ on $\pi_{1}(U_{i})$ for $i=1,2$.
\end{prop}

\subsection{Functorialities and comparison with the case of vector bundles}

In this section we collect some results on the functorial behaviour of the constructed parallel transport and show that our construction is compatible with that in the case of vector bundles in \cite{De-We1}. If we fix a model $\mathfrak{X}$ of $X$, then the functor
$\mathscr{B}_{\mathfrak{X}_{\mathfrak{o}},D}(G)\stackrel{\rho^{\mathfrak{X}}}{\longrightarrow}\Rep_{\pi_{1}(U)}(G(\mathfrak{o}))$ commutes with Galois conjugation, changing the curve or the model, extensions of the structure groups and the action of $\Gal(\overline{\mathbb{Q}}_{p}/K)$ on the given objects in the case that all of them are already defined over a finite extension $K$ of $\mathbb{Q}_{p}$. More details on this functorial behaviour and proofs of it are stated in \cite{H}.

Moreover, one can relate $G$-torsors to vector bundles as follows: Let $\Gamma$ be a free $\mathfrak{o}$-module of finite rank and $\bf{\Gamma}$ the associated affine space over $\spec\,\mathfrak{o}$. Then we define an algebraic representation of $G$ on $\Gamma$ as a homomorphism of $\mathfrak{o}$-group schemes $G\to\GL(\bf{\Gamma})$. In the following we write $\bf{\Gamma}$ for this representation and define $\Rep_{G}(\mathfrak{o})$ and $\Rep_{G(\mathfrak{o})}(\mathfrak{o})$ as the categories of algebraic (resp. continous) representations of $G$ (resp. $G(\mathfrak{o})$) on free $\mathfrak{o}$-modules of finite rank.

For any $G$-torsor $P$ on a $\mathfrak{o}$-scheme $\xi$, the fppf-sheaf \mbox{$E=P\times^{G}\bf{\Gamma}$} is a vector bundle on $\xi$. This construction is functorial in $P$ and the representation $\bf{\Gamma}$; for a morphism $f\colon\widetilde{\xi}\to\xi$ it holds $f^{*}(P\times^{G}{\bf{\Gamma}})=f^{*}P\times^{G}\bf{\Gamma}$. This implies that $\_\times^{G}\bf{\Gamma}$ restricts to a functor $\_\times^{G}{\bf{\Gamma}}\colon\normalfont\mathscr{B}_{\mathfrak{X}_{\mathfrak{o}},D}(G)\to\mathscr{B}_{\mathfrak{X}_{\mathfrak{o}},D}$ such that we obtain a bifunctor \mbox{$\_\times^{G}\_\colon\mathscr{B}_{\mathfrak{X}_{\mathfrak{o}},D}(G)\times\Rep_{G}(\mathfrak{o})\to\mathscr{B}_{\mathfrak{X}_{\mathfrak{o}},D}$}. (As in \cite{De-We1}, $\mathscr{B}_{\mathfrak{X}_{\mathfrak{o}},D}$ is here the category of all vector bundles $E$ on $\mathfrak{X}_{\mathfrak{o}}$ with the property that for all $n\in\N$ exists an element $\pi\in\Sh_{\mathfrak{X},D}$ such that $\pi_{n}^{*}E$ is trivial.)

If $\mathcal{Z}\colon\pi_{1}(U)\to\mathcal{P}(G(\mathfrak{o}))$ is a continous functor and $\Gamma$ a continous representation of $\Gamma$ on $G(\mathfrak{o})$, we define the functor \mbox{$\mathcal{Z}\times^{G(\mathfrak{o})}\Gamma\colon\pi_{1}(U)\to\Mod_{\mathfrak{o}}$} as follows. For $x\in U(\mathbb{C}_{p})$ one sets  $(\mathcal{Z}\times^{G(\mathfrak{o})}\Gamma)(x)=\mathcal{Z}(x)\times^{G(\mathfrak{o})}\Gamma$. For an \'etale path $\gamma$ from $x$ to $x'$ we define $(\mathcal{Z}\times^{G(\mathfrak{o})}\Gamma)(\gamma)$ as the $\mathfrak{o}$-linear map from $\mathcal{Z}(x)\times^{G(\mathfrak{o})}\Gamma$ to $\mathcal{Z}(x')\times^{G(\mathfrak{o})}\Gamma$ which is induced by the continous $G(\mathfrak{o})$-linear map \mbox{$\mathcal{Z}(\gamma)\colon\mathcal{Z}(x)\to\mathcal{Z}(x')$}.

Therefore we may define a bifunctor \[\_\times^{G(\mathfrak{o})}\_\colon\Rep_{\pi_{1}(U)}(G(\mathfrak{o}))\times\Rep_{G(\mathfrak{o})}(\mathfrak{o})\to\Rep_{\pi_{1}(U)}(\mathfrak{o})\] as follows. A pair $(\mathcal{Z},\Gamma)$ of a continous functor $\mathcal{Z}\colon\pi_{1}(U)\to\mathcal{P}(G(\mathfrak{o}))$ and a continous representation of $\Gamma$ on $G(\mathfrak{o})$ is mapped onto the functor \mbox{$\mathcal{Z}\times^{G(\mathfrak{o})}\Gamma\colon\pi_{1}(U)\to\Mod_{\mathfrak{o}}$}.
Morphisms in $\Rep_{\pi_{1}(U)}(G(\mathfrak{o}))$ resp. $\Rep_{G(\mathfrak{o})}(\mathfrak{o})$ are mapped in the obvious way. Like in \cite{De-We1}, $\Rep_{\pi_{1}(U)}(\mathfrak{o})$ is here the $\mathfrak{o}$-linear category of continous functors from $\pi_{1}(U)$ into the category of free $\mathfrak{o}$-modules of finite rank.

\begin{prop}\label{darstellung}
In the above situation, the following diagram commutes:
\begin{equation*}
        \xymatrix@=3em{%
        \mathscr{B}_{\mathfrak{X}_{\mathfrak{o}},D}(G)\times\Rep_{G}(\mathfrak{o}) \ar[rrr]^{\_\times^{G}\_}\ar[d]_{\rho\times\text{nat}} &&& \mathscr{B}_{\mathfrak{X}_{\mathfrak{o}},D} \ar[d]^{\rho}\\
        \Rep_{\pi_{1}(U)}(G(\mathfrak{o}))\times\Rep_{G(\mathfrak{o})}(\mathfrak{o})   \ar[rrr]^{\_\times^{G(\mathfrak{o})}\_} & &&\Rep_{\pi_{1}(U)}(\mathfrak{o}).}
           \end{equation*}
           In particular, $\rho_{P\times^{G}\bf{\Gamma}}=\rho_{P}\times^{G(\mathfrak{o})}\Gamma$ holds as functors from $\pi_{1}(U)$ to $\Mod_{\mathfrak{o}}$ for all $P\in\mathscr{B}_{\mathfrak{X}_{\mathfrak{o}},D}(G)$ and $\Gamma\in\Rep_{G}(\mathfrak{o})$.
\end{prop}
 
\begin{rem}
For any given $G$-torsor $P$ in $\mathscr{B}_{\mathfrak{X}_{\mathfrak{o}},D}(G)$ we obtain through the bifunctor $\_\times^{G}\_$ a vector bundle in $\mathscr{B}_{\mathfrak{X}_{\mathfrak{o}},D}$ for every algebraic representation $\bf{\Gamma}$ of $G$. Thus Proposition \ref{darstellung} shows how to compute the parallel transport along \'etale paths between the fibres of the vector bundle by the parallel transport of the torsor. 
\end{rem}

We explain now shortly how we obtain a principal $\GL_{n}$-bundle in $\mathscr{B}_{\mathfrak{X}_{\mathfrak{o}},D}(\GL_{n})$ from a vector bundle of rank $n$ in $\mathscr{B}_{\mathfrak{X}_{\mathfrak{o}},D}$. If $E$ is a vector bundle of rank $n$ on a scheme $\xi$, then the fppf-sheaf $\Iso_{\xi}(E,\mathbb{A}^{n}_{\xi})$ of local isomorphisms from $E$ to the trivial vector bundle $\mathbb{A}_{\xi}^{n}$ on $\xi$ is a right sheaf torsor under $\GL_{n,\xi}$. Since $\GL_{n,\xi}$ is affine over $\xi$, this sheaf torsor is representable by a $\GL_{n}$-torsor ${\bf{P}}(E)$, the so called \textit{frame bundle} of $E$. If we map any morphism $E\to E'$ to the induced morphism $\Iso_{\xi}(E,\mathbb{A}^{n}_{\xi})\to\Iso_{\xi}(E',\mathbb{A}^{n}_{\xi})$, then we obtain a well-defined functor $\bf{P}$ from the category of vector bundles of rank $n$ on $\xi$ into the category of $\GL_{n}$-torsors on $\xi$. By the construction of the functor holds $f^{*}{\bf{P}}(E)={\bf{P}}(f^{*}E)$ functorially for every morphism $f\colon\widetilde{\xi}\to\xi$. This implies that $\bf{P}$ restricts to a functor from the full subcategory $\mathscr{B}_{\mathfrak{X}_{\mathfrak{o}},D}^{(n)}$ of vector bundles of rank $n$ in $\mathscr{B}_{\mathfrak{X}_{\mathfrak{o}},D}$ into the category $\mathscr{B}_{\mathfrak{X}_{\mathfrak{o}},D}(\GL_{n})$.

On the level of representations, we define a functor $\bf{P}$ from the full subcategory $\Rep_{\pi_{1}(U)}^{n}(\mathfrak{o})$ of $\Rep_{\pi_{1}(U)}(\mathfrak{o})$ which contains of the representations of rank $n$ into $\Rep_{\pi_{1}(U)}(\GL_{n}(\mathfrak{o}))$ as follows:
For an object $\bf{\Gamma}$ of $\Rep_{\pi_{1}(U)}^{n}(\mathfrak{o})$, i.\,e. a continous functor ${\bf{\Gamma}}\colon\pi_{1}(U)\to\Mod_{\mathfrak{o}}^{(n)}$, one defines a continous functor ${\bf{P}}({\bf{\Gamma}})\colon\pi_{1}(U)\to\mathcal{P}(G(\mathfrak{o}))$ by \mbox{${\bf{P}}(\Gamma)(x)=\Iso_{\mathfrak{o}}({\bf}{\Gamma}(x),\mathfrak{o}^{n})$} for every point \mbox{$x\in U(\mathbb{C}_{p})$} where $\Iso_{\mathfrak{o}}(\Gamma(x),\mathfrak{o}^{n})$ is equipped with the naturally induced topology and the right action of $\GL_{n}(\mathfrak{o})$. For an \'etale path $\gamma$ from $x$ to $x'$ one sets  \[{\bf{P}}({\bf{\Gamma}})(\gamma)=({\bf{\Gamma}}(\gamma)^{-1})^{*}\colon\Iso_{\mathfrak{o}}({\bf{\Gamma}}(x),\mathfrak{o}^{n})\to\Iso_{\mathfrak{o}}({\bf{\Gamma}}(x'),\mathfrak{o}^{n}),\varphi\mapsto\varphi\circ\Gamma(\gamma)^{-1}\] such that we obtain altogether a well-defined functor $\bf{P}(\bf{\Gamma})$.

This definition of the functor $\bf{P}$ on the level of objects extends in a natural way to morphisms $\alpha\colon\bf{\Gamma}\to\bf{\Gamma'}$ in $\Rep_{\pi_{1}(U)}^{(n)}(\mathfrak{o})$ by defining the map from $\Iso_{\mathfrak{o}}({\bf{\Gamma}}(x),\mathfrak{o}^{n})$ to $\Iso_{\mathfrak{o}}({\bf{\Gamma'}}(x),\mathfrak{o}^{n})$ which is induced by $\alpha$ for all \mbox{$x\in U(\mathbb{C}_{p})$} as ${\bf{P}}(\alpha)(x)$ which is a natural transformation $\bf{P}(\alpha)\colon\bf{P}(\bf{\Gamma})\to\bf{P}(\bf{\Gamma'})$. 

\begin{prop}\label{pblvbl}
In the given situation, the diagram \begin{equation*}
        \xymatrix@=3em{%
        \mathscr{B}_{\mathfrak{X}_{\mathfrak{o}},D}^{(n)} \ar[r]^{\rho^{\mathfrak{X}}} \ar[d]_{\bf{P}} & \Rep_{\pi_{1}(U)}^{(n)}(\mathfrak{o}) \ar[d]^{\bf{P}}\\
        \mathscr{B}_{\mathfrak{X}_{\mathfrak{o}},D}(\GL_{n}) \ar[r]^{\rho^{\mathfrak{X}}} & \Rep_{\pi_{1}(U)}(\GL_{n}(\mathfrak{o}))}
           \end{equation*} of categories and functors is commutative up to canonical isomorphisms of functors.
           In particular holds $\rho_{{\bf{P}}(E)}={\bf{P}}(\rho_{E})$ for every vector bundle $E$ of rank $n$ in $\mathscr{B}_{\mathfrak{X}_{\mathfrak{o}},D}$. For all $x\in U(\mathbb{C}_{p})$, the representations from $\pi_{1}(U)$ on ${\bf{P}}(E)_{x}=\Iso_{\mathfrak{o}}(E_{x_{\mathfrak{o}}},\mathfrak{o}^{n})$ through $\rho_{{\bf{P}}(E)}$ and through $\gamma\mapsto(\rho_{E}(\gamma)^{-1})^{*}$ are thus equal.
\end{prop}

\begin{rem}
 The Proposition \ref{pblvbl} shows how to compute the isomorphisms of parallel transport for a vector bundle in the category $\mathscr{B}_{\mathfrak{X}_{\mathfrak{o}},D}$ which were constructed by Deninger and Werner (\cite{De-We1}) by those for its frame bundle defined in this article. It follows moreover from the compatiblity of the construction with extensions of structure groups that is suffices to calcute the isomorphisms of parallel tranport for some reduction of the structure group of the frame bundle or equivalently for some subbundle of the given vector bundle.
\end{rem}

Because of the definition of the category $\mathscr{B}_{X_{\mathbb{C}_{p}}, D}(G)$ through the category $\mathscr{B}_{\mathfrak{X}_{\mathfrak{o}},D}(G)$ and by the fact that the functor $\mathscr{B}_{X_{\mathbb{C}_{p}}, D}(G)\to\Rep_{\pi_{1}(U)}(G(\mathbb{C}_{p}))$ was constructed in Section 2.2 in the same way based on the functor $\mathscr{B}_{X_{\mathfrak{o}}, D}(G)\to\Rep_{\pi_{1}(U)}(G(\mathfrak{o}))$, all properties of the category $\mathscr{B}_{\mathfrak{X}_{\mathfrak{o}},D}(G)$ and the functors in the situation over $\mathfrak{o}$ carry over to those over $\mathbb{C}_{p}$. In particular, the functor
\mbox{$\mathscr{B}_{X_{\mathbb{C}_{p}},D}(G)\rightarrow\Rep_{\pi_{1}(U)}(G(\mathfrak{o}))$} commutes with Galois conjugation, changing the curve or the model, extensions of structure groups and the action of $\Gal(\overline{\mathbb{Q}}_{p}/K)$ on the given objects if all of them are already defined over a finite extension $K$ of $\mathbb{Q}_{p}$. Again, more details and proofs of this functorial behaviour are stated in \cite{H}.

Next we state that for any algebraic representation ${\bf{\Gamma}}$ of $G$ the above defined functor $\_\times^{G}{\bf{\Gamma}}\colon\mathscr{B}_{\mathfrak{X}_{\mathfrak{o}},D}(G)\to\mathscr{B}_{\mathfrak{X}_{\mathfrak{o}},D}$ induces a functor \mbox{$\_\times^{G}{\bf{\Gamma}}\colon\mathscr{B}_{X_{\mathbb{C}_{p}}, D}(G)\to\mathscr{B}_{X_{\mathbb{C}_{p}}, D}$}. In the same way one obtains from the bifunctors  \[\_\times^{G}\_\colon\mathscr{B}_{\mathfrak{X}_{\mathfrak{o}},D}(G)\times\Rep_{G}(\mathfrak{o})\to\mathscr{B}_{\mathfrak{X}_{\mathfrak{o}},D}\] and \[\_\times^{G(\mathfrak{o})}\_\colon\Rep_{\pi_{1}(U)}(G(\mathfrak{o}))\times\Rep_{G(\mathfrak{o})}(\mathfrak{o})\to\Rep_{\pi_{1}(U)}(\mathfrak{o})\] bifunctors \[\_\times^{G}\_\colon\mathscr{B}_{X_{\mathbb{C}_{p}}, D}(G)\times\Rep_{G}(\mathfrak{o})\to\mathscr{B}_{X_{\mathbb{C}_{p}}, D}\] and \[\_\times^{G(\mathbb{C}_{p})}\_\colon\Rep_{\pi_{1}(U)}(G(\mathbb{C}_{p}))\times\Rep_{G(\mathbb{C}_{p})}(\mathbb{C}_{p})\to\Rep_{\pi_{1}(U)}(\mathbb{C}_{p})\] where $\Rep_{G(\mathbb{C}_{p})}(\mathbb{C}_{p})$ is the category of continous representations of $G(\mathbb{C}_{p})$ on finite dimensional $\mathbb{C}_{p}$-vector spaces. With these functors follows from proposition \ref{darstellung}:

\begin{prop}\label{darstellung1}
In the above situation, the following diagram commutes:
\begin{equation*}
        \xymatrix@=3em{%
        \mathscr{B}_{X_{\mathbb{C}_{p}},D}(G)\times\Rep_{G}(\mathfrak{o}) \ar[rrr]^{\_\times^{G}\_}\ar[d]_{\rho\times\text{nat}} &&& \mathscr{B}_{X_{\mathbb{C}},D} \ar[d]^{\rho}\\
        \Rep_{\pi_{1}(U)}(G(\mathbb{C}_{p}))\times\Rep_{G(\mathbb{C}_{p})}(\mathbb{C}_{p})   \ar[rrr]^{\_\times^{G(\mathbb{C}_{p})}\_} & &&\Rep_{\pi_{1}(U)}(\mathbb{C}_{p}).}
           \end{equation*}
           In particular,
    $\rho_{P\times^{G}\bf{\Gamma}}=\rho_{P}\times^{G(\mathbb{C}_{p})}(\Gamma\otimes\mathbb{C}_{p})$ holds for all $P\in\mathscr{B}_{X_{\mathbb{C}_{p}},D}(G)$ and \mbox{$\Gamma\in\Rep_{G}(\mathfrak{o})$} as functors from $\pi_{1}(U)$ to $\mathcal{P}(G(\mathbb{C}_{p}))$.
\end{prop}

Finally, the above defined functors ${\bf{P}}$ extend to functors \mbox{${\bf P}\normalfont\colon\mathscr{B}_{X_{\mathbb{C}_{p}},D}^{(n)}\rightarrow\mathscr{B}_{X_{\mathbb{C}_{p}},D}(\GL_{n})$} and $\bf{P}\colon\Rep_{\pi_{1}(U)}^{(n)}(\mathbb{C}_{p})\longrightarrow\Rep_{\pi_{1}(U)}(\GL_{n}(\mathbb{C}_{p}))$. From Proposition \ref{pblvbl}, it follows: 

\begin{prop}\label{vb1}
In the given situation, the diagram \begin{equation*}
        \xymatrix@=3em{%
        \mathscr{B}_{X_{\mathbb{C}_{P}},D}^{(n)} \ar[r]^{\rho} \ar[d]_{\bf{P}} & \Rep_{\pi_{1}(U)}^{(n)}(\mathbb{C}_{p}) \ar[d]^{\bf{P}}\\
        \mathscr{B}_{X_{\mathbb{C}_{p}},D}(\GL_{n}) \ar[r]^{\rho} & \Rep_{\pi_{1}(U)}(\GL_{n}(\mathbb{C}_{p}))}
           \end{equation*} of categories and functors is commutative up to canonical isomorphisms of functors.
In particular, $\rho_{{\bf{P}}(E)}={\bf{P}}(\rho_{E})$ holds for every vector bundle $E$ of rank $n$ in $\mathscr{B}_{X_{\mathbb{C}_{p}},D}$. Thus the representations from $\pi_{1}(U)$ to $({\bf{P}}(E))_{x}=\Iso_{\mathbb{C}_{p}}(E_{x},\mathbb{C}_{p}^{n})$ through $\rho_{{\bf{P}}(E)}$ and through $\gamma\mapsto(\rho_{E}(\gamma)^{-1})^{*}$ are equal for all $x\in U(\mathbb{C}_{p})$.
\end{prop}

\section{Characterisation of the category $\mathscr{B}_{\mathfrak{X}_{\mathfrak{o}},D}(G)$}

It is difficult to verify directly from Definition \ref{modcat} whether a given $G$-torsor $P$ on $\mathfrak{X}_{\mathfrak{o}}$ is contained in the category $\mathscr{B}_{\mathfrak{X}_{\mathfrak{o}},D}(G)$. In the following we establish therefore another characterisation of the objects of $\mathscr{B}_{\mathfrak{X}_{\mathfrak{o}},D}(G)$. This is analogous to Theorems 16, 17 and 18 in \cite{De-We1}.

\subsection{Reduction to the special fibre}

\begin{thm}\label{thm16}
Let $\mathfrak{X}$ be a model over $\overline{\mathbb{Z}}_{p}$ of a smooth projective curve $X$ over $\overline{\mathbb{Q}}_{p}$. and denote by $k=\overline{\mathbb{F}}_{p}$ the residue field of $\overline{\mathbb{Z}}_{p}$. A $G$-torsor $P$ on $\mathfrak{X}_{\mathfrak{o}}$ is contained in $\mathscr{B}_{\mathfrak{X}_{\mathfrak{o}},D}(G)$ if and only if there is an object $\pi\colon\mathcal{Y}\to\mathfrak{X}$ of the category $\Sh_{\mathfrak{X},D}$ such that $\pi_{k}^{*}P_{k}$ is trivial on \mbox{$\mathcal{Y}_{k}=\mathcal{Y}\otimes_ {\overline{\mathbb{Z}}_{p}}k$} where we set $P_{k}=P\otimes_{\mathfrak{o}}k$.
\end{thm}
\iproof{\begin{itemize}{\item[(i)] The necessity is obvious.
Let thus $P$ be a $G$-torsor on $\mathfrak{X}_{\mathfrak{o}}$ with the property that there is an object \mbox{$\pi\colon\mathcal{Y}\to\mathfrak{X}$} of $\Sh_{\mathfrak{X},D}$ such that $\pi_{k}^{*}P_{k}$ is trivial on \mbox{$\mathcal{Y}_{k}=\mathcal{Y}\otimes_ {\overline{\mathbb{Z}}_{p}}k$}. By Remark \ref{doms} we may assume $\pi\in S_{\mathfrak{X},D}^{\text{ss}}$.
Then it follows from noetherian descent that the family $(\mathfrak{X}, D, G, P_{1},\pi\colon\mathcal{Y}\to\mathfrak{X})$ descends to a family \mbox{$(\mathfrak{X}^{0},D^{0},G^{0}, \overline{P}, \pi_{0}\colon\mathcal{Y}^{0}\to\mathfrak{X}^{0})$} over the ring of integers $\mathfrak{0}_{K}$ of a finite extension $K$ of $\mathbb{Q}_{p}$.
Here $\mathfrak{X}^{0}$ is a model over $\mathfrak{o}_{K}$ of the smooth projective curve \mbox{$X^{0}=\mathfrak{X}^{0}\otimes_{\mathfrak{o}_{K}}K$} over $K$, $G^{0}$ a group scheme which is affine, smooth and of finite presentation over a normal, finitely generated $\mathfrak{o}_{K}$-algebra $A$ and $\overline{P}$ is a $G^{0}_{\mathfrak{X}^{0}_{A}\otimes_{A}\mathfrak{o}_{K}/p\mathfrak{o}_{K}}$-torsor on $\mathfrak{X}^{0}_{A}\otimes_{A}\mathfrak{o}_{K}/p\mathfrak{o}_{K}$ whose restriction to the special fibre \mbox{$\mathfrak{X}^{0}\otimes_{\mathfrak{o}_{K}}\mathfrak{o}_{K}/\mathfrak{p}$} becomes trivial after pullback along $\pi_{0}\otimes\mathfrak{o}_{K}/\mathfrak{p}$. Moreover, $\pi_{0}$ is an element of $S_{\mathfrak{X}_{0},D_{0}}^{\text{ss}}$.
 \item[(ii)] Let $e$ be the ramification index of $K$ over $\mathbb{Q}_{p}$ and let \mbox{$\mathfrak{o}_{\nu/e}:=\mathfrak{o}/\mathfrak{p}^{\nu}\mathfrak{o}=\overline{\mathbb{Z}}_{p}/\mathfrak{p}^{\nu}\overline{\mathbb{Z}}_{p}$}. Note that this is compatible with our earlier notation \mbox{$\mathfrak{o}_{n}=\mathfrak{o}/p^{n}\mathfrak{o}$}.
 
Let $\pi_{\nu/e}, \widetilde{P}_{\nu/e}$ etc. be always the base change with $\mathfrak{o}_{\nu/e}$ where $\widetilde{P}$ is the generic fibre of $P$. As $\pi_{1/e}$ is also the base change of \mbox{$\pi_{0}\otimes_{\mathfrak{o}_{K}}\mathfrak{o}_{K}/\mathfrak{p}$} with $\mathfrak{o}_{1/e}$, it follows that $\pi_{1/e}^{*}P_{1/e}=(\pi_{0}\otimes\mathfrak{o}_{K}/\mathfrak{p}\otimes\mathfrak{o}_{1/e})^{*}(\overline{P}\otimes\mathfrak{o}/\mathfrak{p})$ is trivial on $\mathcal{Y}_{1/e}$.

By induction it therefore suffices to prove the following assertion (*):

Given $\nu\geq 2$ and some element $\pi\colon\mathcal{Y}\to\mathfrak{X}$ of $\Sh_{\mathfrak{X},D}^{\text{ss}}$ such that $\pi^{*}_{(\nu-1)/e}P_{(\nu-1)/e}$ is trivial, there exists an object $\mu\colon\mathcal{Z}\to\mathfrak{X}$ in $\Sh_{\mathfrak{X},D}^{\text{ss}}$ such that $\mu^{*}_{\nu/e}P_{\nu/e}$ is a trivial $G$-torsor on $\mathcal{Z}_{\nu/e}$. 

Let $\nu$ and \mbox{$\pi\colon\mathcal{Y}\rightarrow\mathfrak{X}$} be as in $(*)$ and \mbox{$i\colon\mathcal{Y}_{(\nu-1)/e}\lhook\rightarrow\mathcal{Y}_{\nu/e}$} the natural closed immersion.

If $\mathcal{Y}_{\nu/e}$ is affine, then the canonical application \[H^{1}(\mathcal{Y}_{\nu/e},G_{\mathcal{Y}_{\nu/e}})\to H^{1}(\mathcal{Y}_{(\nu-1)/e},G_{\mathcal{Y}_{(\nu-1)/e}})\] is by \cite[Corollaire VII.1.3.2]{Gi} a bijection as for fppf topology as for \'etale topology and thus $\pi^{*}_{\nu/e}P_{\nu/e}$ is also trivial such that we may choose $\mu=\pi$.

In the following we assume therefore that $\mathcal{Y}_{\nu/e}$ is not affine. Let $\mathcal{I}$ be the sheaf of ideals which defines the immersion $i$ and satisfies $\mathcal{I}^{2}=0$.
Then we have by \cite[Lemme VII.1.3.5]{Gi} the following exact sequence of sheaves on $(\mathcal{Y}_{\nu/e})_{\acute{e}t}$: \begin{equation*}
        \xymatrix@=1em{%
         1 \ar[r] & i_{*}(\mathscr{L}_{(\nu-1)/e}\otimes_{\mathcal{O}_{\mathcal{Y}_{(\nu-1)/e}}}\mathcal{I}) \ar[rrrr]^-{h} & & & & G \ar[rrr]^-{adj}& & & i_{*}G_{\mathcal{Y}_{(\nu-1)/e}} \ar[r] & 1,}
           \end{equation*}
 where we set $\mathscr{L}_{(\nu-1)/e}:=Lie(G_{\mathcal{Y}_{(\nu-1)/e}}/\mathcal{Y}_{(\nu-1)/e})$.        Note that $i$ is an universal homoeomorphism such that 
one may identify $\mathcal{I}$ naturally with an $\mathcal{O}_{\mathcal{Y}_{(\nu-1)/e}}$-module. Moreover, $\adj$ denotes the composition of the adjunction morphism \mbox{$G_{\mathcal{Y}_{\nu/e}}\rightarrow i_{*}i^{*}G_{\mathcal{Y}_{\nu/e}}$} with the morphism $i_{*}i^{*}G_{\mathcal{Y}_{\nu/e}}\rightarrow i_{*}G_{\mathcal{Y}_{(\nu-1)/e}}$ induced by the canonical morphism \mbox{$i^{*}G_{\mathcal{Y}_{\nu/e}}\rightarrow G_{\mathcal{Y}_{(\nu-1)/e}}$}. Next, the morphism $h$ is induced by the bijection \mbox{$i_{*}(\mathscr{L}_{(\nu-1)/e}\otimes\mathcal{I})(S)\to\ker(G_{\mathcal{Y}_{\nu/e}}\to i_{*}G_{\mathcal{Y}_{(\nu-1)/e}})(S),\alpha\to e\alpha$}. Here $e$ is the unit section of $G_{\mathcal{Y}_{(\nu-1)/e}}$ and $S\to\mathcal{Y}_{\nu/e}$ is \'etale.

From this exact sequence we deduce the exact cohomology sequence \[H^{1}_{\acute{e}t}(\mathcal{Y}_{\nu/e},i_{*}(\mathscr{L}_{(\nu-1)/e}\otimes_{\mathcal{O}_{\mathcal{Y}_{(\nu-1)/e}}}\mathcal{I}))\xrightarrow{h} H^{1}_{\acute{e}t}(\mathcal{Y}_{\nu/e},G_{\mathcal{Y}_{\nu/e}})\to H^{1}_{\acute{e}t}(\mathcal{Y}_{\nu/e},i_{*}G_{\mathcal{Y}_{(\nu-1)/e}}).\]
By \cite[Corollaire VIII.5.6]{SGA 4} is now $R^{1}i_{*}G_{\mathcal{Y}_{(\nu-1)/e}}=0$ because $i$ is a closed immersion and in particular integral. This implies by \cite[Proposition V.3.1.3]{Gi} that $H^{1}_{\acute{e}t}(\mathcal{Y}_{\nu/e},i_{*}G_{\mathcal{Y}_{(\nu-1)/e}})$ is isomorphic to $H^{1}_{\acute{e}t}(\mathcal{Y}_{(\nu-1)/e},G_{\mathcal{Y}_{(\nu-1)/e}})$ such that the sequence  \[H^{1}_{\acute{e}t}(\mathcal{Y}_{\nu/e},i_{*}(\mathscr{L}_{(\nu-1)/e}\otimes_{\mathcal{O}_{\mathcal{Y}_{(\nu-1)/e}}}\mathcal{I}))\xrightarrow{h} H^{1}_{\acute{e}t}(\mathcal{Y}_{\nu/e},G_{\mathcal{Y}_{\nu/e}})\xrightarrow{i^{*}} H^{1}_{\acute{e}t}(\mathcal{Y}_{(\nu-1)/e},G_{\mathcal{Y}_{(\nu-1)/e}})\] is therefore exact.
As above, $R^{1}i_{*}(\mathscr{L}_{(\nu-1)/e}\otimes_{\mathcal{O}_{\mathcal{Y}_{(\nu-1)/e}}}\mathcal{I})=0$ holds by \cite[Corollaire VIII.5.6]{SGA 4}. Thus we have by \cite[Proposition V.3.1.3]{Gi} an isomorphism \begin{equation*}
 \xymatrix@=3em{%
H^{1}_{\acute{e}t}(\mathcal{Y}_{(\nu-1)/e},\mathscr{L}_{(\nu-1)/e}\otimes_{\mathcal{O}_{\mathcal{Y}_{(\nu-1)/e}}}\mathcal{I})\ar[r]^{j}_{\sim} & H^{1}_{\acute{e}t}(\mathcal{Y}_{\nu/e},i_{*}(\mathscr{L}_{(\nu-1)/e}\otimes_{\mathcal{O}_{\mathcal{Y}_{(\nu-1)/e}}}\mathcal{I}))}\end{equation*} such that the sequence  \[H^{1}_{\acute{e}t}(\mathcal{Y}_{(\nu-1)/e},\mathscr{L}_{(\nu-1)/e}\otimes_{\mathcal{O}_{\mathcal{Y}_{(\nu-1)/e}}}\mathcal{I})\xrightarrow{h\circ j} H^{1}_{\acute{e}t}(\mathcal{Y}_{\nu/e},G_{\mathcal{Y}_{\nu/e}})\xrightarrow{i^{*}} H^{1}_{\acute{e}t}(\mathcal{Y}_{(\nu-1)/e},G_{\mathcal{Y}_{(\nu-1)/e}})\] is therefore exact.

If $\Omega$ is the class of $\pi_{\nu/e}^{*}P_{\nu/e}$ in $H^{1}_{\acute{e}t}(\mathcal{Y}_{\nu/e},G_{\mathcal{Y}_{\nu/e}})$, then $i^{*}$ maps $\Omega$ onto the class of $i^{*}\pi_{\nu/e}^{*}P_{\nu/e}=\pi^{*}_{(\nu-1)/e}P_{(\nu-1)/e}$, i.\,e. the trivial class in $H^{1}_{\acute{e}t}(\mathcal{Y}_{(\nu-1)/e},G_{\mathcal{Y}_{(\nu-1)/e}})$.

As the above sequence is exact, there exists a class $A$ in \mbox{$H^{1}_{\acute{e}t}(\mathcal{Y}_{(\nu-1)/e},\mathscr{L}_{(\nu-1)/e}\otimes_{\mathcal{Y}_{(\nu-1)/e}}\mathcal{I})$} such that $(h\circ j)(A)=\Omega$. The canonical $\mathcal{O}_{\mathcal{Y}_{(\nu-1)/e}}$-module structure of $\mathscr{L}_{(\nu-1)/e}$ and $\mathcal{I}$ induces moreover the following exact sequence of sheaves on $(\mathcal{Y}_{(\nu-1)/e})_{\acute{e}t}$: 
\[0\to\ker(\alpha)\to\mathcal{I}\xrightarrow{\alpha}\mathscr{L}_{(\nu-1)/e}\otimes_{\mathcal{O}_{\mathcal{Y}_{(\nu-1)/e}}}\mathcal{I}\to0.\] Here $\alpha$ is the canonical morphism. This delivers a surjection \begin{equation*}\xymatrix@=2em{%
H^{1}_{\acute{e}t}(\mathcal{Y}_{(\nu-1)/e},\mathcal{I})\ar@{>>}[rrr]^-{\alpha} & & & H^{1}_{\acute{e}t}(\mathcal{Y}_{(\nu-1)/e},\mathscr{L}_{(\nu-1)/e}\otimes_{\mathcal{Y}_{(\nu-1)/e}}\mathcal{I})}\end{equation*} because $\mathcal{Y}_{(\nu-1)/e}$ is onedimensional and  $H^{2}_{\acute{e}t}(\mathcal{Y}_{(\nu-1)/e},\mathscr{F})\cong H^{2}(\mathcal{Y}_{(\nu-1)/e},\mathscr{F}_{Zar})$ holds for every \'etale sheaf $\mathscr{F}$ on $\mathcal{Y}_{(\nu-1)/e}$ by \cite[Proposition III.3.3]{M} where $\mathscr{F}_{Zar}$ is the restriction to the Zariski topology.

The $\mathcal{O}_{\mathcal{Y}_{(\nu-1)/e}}$-module structure of $\mathcal{I}$ gives us moreover a canonical exact sequence \[0\longrightarrow\ker g\longrightarrow\mathcal{O}_{\mathcal{Y}_{(\nu-1)/e}}\stackrel{g}{\longrightarrow}\mathcal{I}\longrightarrow 0\] of sheaves on $(\mathcal{Y}_{(\nu-1)/e})_{\acute{e}t}$. The same argument as above delivers thus a surjection $H^{1}_{\acute{e}t}(\mathcal{Y}_{(\nu-1)/e},\mathcal{O}_{\mathcal{Y}_{(\nu-1)/e}})\stackrel{g}{\twoheadrightarrow} H^{1}_{\acute{e}t}(\mathcal{Y}_{(\nu-1)/e},\mathcal{I})$. 

If next the genus of the generic fibre $Y_{K}$ of the model $\mathcal{Y}^{0}$ is zero, then it follows in the same way as in \cite{De-We1} that $H^{1}_{Zar}(\mathcal{Y}_{(\nu-1)/e},\mathcal{O}_{\mathcal{Y}_{(\nu-1)/e}})=0$ holds: As we may assume that $Y_{K}$ has a rational point, we have \mbox{$Y_{K}=\mathbb{P}^{1}_{K}$} such that \mbox{$\chi(Y_{K},\mathcal{O}_{Y_{K}})=\chi(\mathcal{Y}_{\kappa},\mathcal{O}_{Y_{\kappa}})=1$} where $\mathcal{Y}_{\kappa}$ is the special fibre of $\mathcal{Y}^{0}$. As moreover  $(\lambda_{0})_{*}\mathcal{O}_{\mathcal{Y}_{0}}=\mathcal{O}_{\spec\,\mathfrak{o}_{K}}$ where $\lambda_{0}$ denotes the structure morphism of $\mathcal{Y}_{0}$, this implies  $H^{1}(\mathcal{Y}_{\kappa},\mathcal{O}_{\mathcal{Y}_{\kappa}})=0$. By \cite[p. 53, Corollary 3]{Mf} it holds therefore $H^{1}(\mathcal{Y}_{\nu'/e},\mathcal{O}_{\mathcal{Y}_{\nu'/e}})=0$ for all $\nu'\in\N$. This implies $H^{1}_{Zar}(\mathcal{Y}_{(\nu-1)/e},\mathcal{O}_{\mathcal{Y}_{(\nu-1)/e}})=0$ and by \cite[Proposition I.2.5.]{FK} $H^{1}_{\acute{e}t}(\mathcal{Y}_{(\nu-1)/e},\mathcal{O}_{\mathcal{Y}_{(\nu-1)/e}})=0$ either.
But then $A=0$ holds and thus $\Omega=0$ either such that $\pi_{\nu/e}^{*}P_{\nu/e}$ is  trivial. In the case that the genus of $Y_{K}$ is zero we may therefore choose $\mu=\pi$.

Let us therefore now assume that the genus of $Y_{K}$ is different from zero. Then it follows from the proof of \cite[Theorem 11]{De-We1} that there is a morphism $\rho\colon(\mu\colon\mathcal{Z}\to\mathfrak{X})\to(\pi\colon\mathcal{Y}\to\mathfrak{X})$ in $\Sh_{\mathfrak{X},D}$ such that the induced map $\rho^{*}\colon H^{1}(\mathcal{Y},\mathcal{O}_{\mathcal{Y}})\to H^{1}(\mathcal{Z},\mathcal{O}_{\mathcal{Z}})$ satisfies the relation \[\rho^{*}(H^{1}(\mathcal{Y},\mathcal{O}_{\mathcal{Y}}))\subset p^{\nu/e}H^{1}(\mathcal{Z},\mathcal{O}_{\mathcal{Z}}).\] By \cite[Proposition I.2.5]{FK} one has the same assertion in the case that we work with cohomology groups for the \'etale cohomology instead of for the Zariski topology. If we work in \'etale cohomolgy, we denote the morphism by $\rho_{\acute{e}t}^{*}$. As in the proof of \cite[Theorem 16]{De-We1} it follows that the induced map \[(\rho_{\acute{e}t})_{(\nu-1)/e}^{*}\colon H^{1}_{\acute{e}t}(\mathcal{Y}_{(\nu-1)/e},\mathcal{O}_{\mathcal{Y}_{(\nu-1)/e}})\to H^{1}_{\acute{e}t}(\mathcal{Z}_{(\nu-1)/e},\mathcal{O}_{\mathcal{Z}_{(\nu-1)/e}})\] is trivial. From the commutativity of the diagram \begin{equation*}
 \xymatrix@=3em{%
H^{1}_{\acute{e}t}(\mathcal{Y}_{(\nu-1)/e},\mathcal{O}_{\mathcal{Y}_{(\nu-1)/e}})\ar[rrr]^{h\circ j\circ \alpha\circ g} \ar[d]^{(\rho^{*}_{\acute{e}t})_{(\nu-1)/e}=0} & && H^{1}_{\acute{e}t}(\mathcal{Y}_{\nu/e},G_{\mathcal{Y}_{\nu/e}}) \ar[d]^{(\rho^{*}_{\acute{e}t})_{\nu/e}}\\
H^{1}_{\acute{e}t}(\mathcal{Z}_{(\nu-1)/e},\mathcal{O}_{\mathcal{Z}_{(\nu-1)/e}})\ar[rrr]^{h\circ j\circ \alpha\circ g} & & & H^{1}_{\acute{e}t}(\mathcal{Z}_{\nu/e},G_{\mathcal{Y}_{\nu/e}})}\end{equation*} follows then that $(\rho_{\acute{e}t})_{\nu/e}^{*}\Omega$ is trivial because, as already shown, $\Omega$ possedes a preimage in $H^{1}_{\acute{e}t}(\mathcal{Y}_{(\nu-1)/e},\mathcal{O}_{\mathcal{Y}_{(\nu-1)/e}})$ under $h\circ j\circ\alpha\circ g$.

Therefore $\mu_{\nu/e}^{*}P_{\nu/e}=(\rho_{\acute{e}t})_{\nu/e}^{*}(\pi_{\nu/e}^{*}P_{\nu/e})$ is also a trivial torsor on $\mathcal{Z}_{\nu/e}$ for the \'etale topology and for the fppf-topology what was to show.}\endproof\end{itemize}

\subsection{Trivializibility of \'etale $G$-torsors}

In order to find a characterisation of the category $\mathscr{B}_{\mathfrak{X}_{\mathfrak{o}},D}(G)$ without any use of the covering categories $\Sh_{\mathfrak{X},D}, S_{\mathfrak{X},D}^{\text{good}}$ resp. $S^{\text{ss}}_{\mathfrak{X},D}$ we need the following definition:

\begin{defn}
Let $G$ be a connected algebraic group over a finite field $\mathbb{F}_{q}$ and $X$ be an $\mathbb{F}_{q}$-scheme. Then we call an \'etale $G$-torsor $P$ trivialisable, if there exists a finite surjective morphism \mbox{$f\colon Y\to X$} such that $f^{*}P$ is a trivial $G$-torsor on $Y$.
\end{defn}

We recall now the following resultat by Deligne cited in \cite{Las} which gives us a relation between trivialisable \'etale $G$-torsors and those that are isomorphic to its pullback under the absolute Frobenius:

\begin{prop}[see \protect{\cite[Proof of Lemma 3.3.]{Las}}]\label{Laszlo}
Let $G$ be a smooth connec\-ted algebraic group over $\mathbb{F}_{q}$ and $X$ be a scheme over $\mathbb{F}_{q}$ and $F$ the absolute Frobenius. Here it is $q=p^{r}$ and our notation does not distinguish in order to simplify her between the absolute Frobenius of $X$ and those of $G$. Then the embedding \mbox{$G(\mathbb{F}_{q})= G^{F}\into G$} induces a functor from the category of $G(\mathbb{F}_{q})$-torsors on $X$ with respect to the \'etale topology into the category of \'etale $G$-torsors $P$ on $X$ which are isomorphic to $F^{*}P$ as $G$-torsors. This functor is an equivalence of categories.  
\end{prop}

\begin{cor}\label{pullback2}
Let $G$ be a connected reductive group scheme of finite presentation over $\mathbb{F}_{q}$, $q=p^{r}$, $X$ be a scheme over $\mathbb{F}_{q}$ and $F$ the absolute Frobenius of $X$. Then every \'etale $G$-torsor $P$ on $X$, which is isomorphic to its pullback $F^{*}P$ as an \'etale $G$-torsor, is trivialisable by a finite \'etale morphism $f\colon Y\to X$. 
\end{cor}
\begin{proof}{ Let $P$ be an \'etale $G$-torsor on $X$ which is as an \'etale $G$-torsor isomorphic to its pullback $F^{*}P$ under the Frobenius. Then there exists following Proposition \ref{Laszlo} a $G^{F}$-torsor $Q$ such that $P=Q\wedge^{G^{F}}G$ holds. As $G^{F}=G(\mathbb{F}_{q})$ is smooth and finite, $Q$ is finite over $X$ by \cite[Proposition III.4.2]{M} such that $f\colon Q\to X$ is a finite \'etale covering of $X$. One should notice here that under the given assumptions every $G$-torsor with respect to the fppf-topology is also a torsor with respect to the \'etale topology and vice versa. It is now $f^{*}Q=Q\times_{X}Q$ trivial as an \'etale $G^{F}\times_{X}Q$-torsor on $Q$ such that \mbox{$f^{*}P=P\times_{X}Q=(Q\times_{X}Q)\wedge^{G^{F}}G$} is also trivial as an \'etale $G$-torsor.}\end{proof}

\begin{cor}\label{potenz}
Let $G$ be a connected reductive group scheme of finite presentation over $\mathbb{F}_{q}$, $q=p^{r}$, $X$ be a scheme over $\mathbb{F}_{q}$ and let $P$ be an \'etale $G$-torsor on $X$ with the property that there is an $s\in\N$ such that $(F^{t})^{*}P\cong(F^{s})^{*}P$ holds for all $t\geq s$ as \'etale $G$-torsors where $F$ denotes the absolute Frobenius of $X$. Then $P$ is trivialisable.
\end{cor} 

\subsection{Semistability of principal bundles}

As in the case of vector bundles these results will enable us a description of the category $\mathscr{B}_{\mathfrak{X}_{\mathfrak{o},D}}(G)$ using the notion of \emph{semistability}. Therefore we recall the following construction by Serre (\cite{Se}):

Given a smooth reductive group $G$ over a field $k$ of arbitrary characteristic,  a smooth projective curve $X$ over $k$, a $G$-torsor $E$ on $X$ and a quasi-projective scheme $F$ over $k$ on which $G$ acts from the left, an action of $G$ on $E\times_{\spec\,k}F$ is defined by $g(e,f):=(e(g\otimes 1),g^{-1}f)$ for all \mbox{$g\in G$}, $e\in E$ and $f\in F$. The quotient of $E\times_{\spec\,k}F$ modulo this action of $G$ is representable by an unique scheme $E(F)$ such that the morphism \mbox{$E\times_{\spec\,k}F\to E(F)$} gives the structure of an \'etale $G$-torsor over $E(F)$ to the fibre product $E\times_{\spec\,k}F$. It is clear that every morphism $\alpha\colon F_{1}\to F_{2}$ of quasi-projective $k$-schemes with an  action from $G$ from the left which is compatible with the action of $G$ on $F_{1}$ resp. $F_{2}$ induces a morphism $\alpha_{E}\colon E(F_{1})\to E(F_{2})$. If $\sigma\colon G\to H$ is a morphism of connected reductive algebraic groups over $k$, the above defined scheme $E_{\rho}(H)$ is an $H$-torsor. If moreover $H$ is a normal closed subgroup of $G$, then it follows from Chevalley's theorem that $G/H$ is a quasi-projective variety. If we equip $G/H$ with the natural action of $G$ from the left, then there exists therefore the above defined scheme $E(G/H)$.

Using this construction by Serre we can now define extensions and reductions of the structure group:

\begin{defn}
If $\rho\colon G\to H$ is a morphism of connected reductive group schemes of finite presentation over a field $k$ and $X$ a smooth projective curve over $k$, then $\rho$ induces a map $H^{1}(X,G)\to H^{1}(X,H)$ by mapping each $G$-torsor $E$ onto the $H$-torsor $E_{\rho}(H)$. One says that one obtains $E_{\rho}(H)$ by extension of the strucure group of $E$ to $H$. If there is on the other hand a $H$-torsor $E'$ given together with a $G$-torsor $E$ and an isomorphism of $H$-torsors \mbox{$\Phi\colon E_{\rho}(H)\to E'$}, then we call $(E,\Phi)$ a reduction of the structure group of $E'$ to $G$.
\end{defn}

We set $E/P:=E(G/P)$ and denote by $T_{E/P}$ the relative tangent bundle with respect to the projection \mbox{$E/P\to X$.} If $\sigma\colon X\to E/P$ is a section of the canonical projection \mbox{$E/P\to X$}, one finds that $\sigma^{*}T_{E/P}=E(\mathfrak{g}/\mathfrak{p})$. Moreover, the sections $\sigma\colon X\to E/P$ are in 1:1-correspondence to the \'etale $P$-torsors on $X$ which one obtains by reductions of the structure group of $E$ to $P$. Thus we may denote reductions of the structure group of $E$ to $P$ in the following by $(P,\sigma)$.

For the upcoming characterisation of the category $\mathscr{B}_{\mathfrak{X}_{\mathfrak{o}}, D}(G)$ we need the notion of \emph{semistability} of principal bundles. This notion was defined by Ramanathan extending the definition of semistability for vector bundles by Mumford as follows:

\begin{defn}[Ramanathan, \protect{\cite[Definition 1.1]{R}}]\label{semistabil}
Let $G$ be a connected reductive algebraic group over a field $k$ and $X$ a smooth projective curve over $k$. Then a $G$-torsor $E$ on $X$ is said to be semistable, if and only if for every reduction $(P,\sigma)$ of the structure group to a maximal parabolic subgroup $P$, where $\sigma$ is a section $\sigma\colon X\to E/P$, the inequality $\deg(\sigma^{*}T_{E/P})\geq 0$ holds. 
\end{defn}

From the relation $\sigma^{*}T_{E/P}=E(\mathfrak{g}/\mathfrak{p})$ it follows moreover:

\begin{lem}[compare with \protect{\cite[Remark 2.2]{R}}, \protect{\cite[Lemma 2.5]{BH}}]\label{adsemi}
Let $G$ be a connected reductive algebraic group over a field $k$ of arbitrary characteristic and $X$ a smooth projective curve over $k$. Then a $G$-torsor $E$ on $X$ is semistable if the associated vector bundle $E(\mathfrak{g})$ is semistable.
\end{lem}
\begin{proof}{
Let $(P,\sigma)$ be a reduction of the structure group of $E$ to a maximal parabolic subgroup $P$ of $G$. Then the exact sequence $0\to\mathfrak{p}\to\mathfrak{g}\to\mathfrak{g}/\mathfrak{p}\to 0$ 
           gives us by the functoriality of the construction of the associated fibre bundle an exact sequence $0\to E(\mathfrak{p})\to E(\mathfrak{g})\to E(\mathfrak{g}/\mathfrak{p})\to 0$
          of vector bundles on $X$. As $E(\mathfrak{g})$ is semistable as a vector bundle, it is therefore \mbox{$\mu(E(\mathfrak{g}/\mathfrak{p}))\geq\mu(E(\mathfrak{g}))$}, i.\,e. \mbox{$\mu(\sigma^{*}T_{E/P})\geq\mu(E(\mathfrak{g}))$}. As by \cite[Note 4.2]{Be} \mbox{$\deg(E(\mathfrak{g}))=0$} holds such that we have $\mu(E(\mathfrak{g}))=0$, it follows \mbox{$\deg(\sigma^{*}T_{E/P})\geq 0$} and therefore the semistability of $E$. }\end{proof}

\begin{cor}\label{torsortrivial}
Let $G$ be a connected reductive algebraic group over a field $k$ and $X$ a smooth projective curve over $k$. Then the trivial $G$-torsor $G_{X}$ on $X$ is semistable.
\end{cor}
\begin{proof}{ In the given situtation holds \[G_{X}(\mathfrak{g})=(G\times_{\spec\,k}X)\times_{\spec\,k}\mathfrak{g}/(x,eg,g^{-1}f)\cong X\times_{\spec\,k}(G\times_{\spec\,k}\mathfrak{g}/(eg,g^{-1}f))\] with $x\in X$, $e\in G$ and $f\in\mathfrak{g}$ such that $G_{X}(\mathfrak{g})$ is a trivial vector bundle. Following Lemma \ref{adsemi}, $G_{X}$ is therefore semistable.}\end{proof}

\begin{rem}
Let $G, X$ and $k$ as before.
Then the proof of Lemma \ref{adsemi} shows also that a $G$-torsor $E$ on $X$ is semistable if and only if  $\deg(\sigma^{*}E(\mathfrak{p}))\leq 0$ holds for any reduction $(P,\sigma)$ of the structure group to a maximal parabolic subgroup $P$ of $G$.
\end{rem}

As in the case of vector bundles one may check semistability also for $G$-torsors after pullback along a finite covering: 

\begin{lem}[compare \protect{\cite[Lemma 6.8]{BH}}]\label{endlichsemist}
Let $X$ be a smooth projective curve over a field $k$, $G$ a connected reductive algebraic group over $k$, $E$ a $G$-torsor on $X$ and $f\colon Y\to X$ a finite covering. Then $E$ is semistable, if  $f^{*}E$ is it.
\end{lem} 
\begin{proof}{ Let $f^{*}E$ be semistable. Suppose that $E$ is not semistable. By Definition \ref{semistabil} exist a maximal parabolic subgroup $P$ of $G$ and a section $\sigma\colon X\to E/P$ such that $\deg(\sigma^{*}T_{E/P})<0$. But \[\deg((\sigma\otimes\id_{Y})^{*}T_{f^{*}E/P_{Y}})=\deg(f^{*}(\sigma^{*}T_{E/P}))=[K(Y):K(X)]\deg(\sigma^{*}T_{E/P})\] implies then $\deg((\sigma\otimes\id_{Y})^{*}T_{E\times_{X}Y/P_{Y}})<0$ and $f^{*}=E\times_{X}Y$ is therefore not semistable in contradiction to our assumption.}\end{proof}

Furthermore it holds:

\begin{lem}\label{noether}
Let $X_{0}$ be a smooth projective curve over a field $k$, $\overline{k}$ an algebraic closure of $k$, $X:=X_{0}\otimes_{k}\overline{k}$ and $f:X\to X_{0}$ the canonical projection. Moreover, let $G$ be a connected reductive algebraic group over $k$ and $E$ a $G$-torsor on $X_{0}$. Then $E$ is semistable, if $f^{*}E$ is semistable.
\end{lem} 
\begin{proof}{ Let $f^{*}E$ be semistable. Suppose that $E$ would be not semistable. Then there exist by Definition \ref{semistabil} a maximal parabolic subgroup $P$ of $G_{X_{0}}$ and a section \mbox{$\sigma\colon X\to E/P$} such that $\deg(\sigma^{*}T_{E/P})<0$. As we have $f^{*}E=E\otimes_{k}\overline{k}$ and all constructions are compatible with base change, it follows by \cite[Proposition 7.3.7]{Liu}  $\deg((\sigma\otimes_{k}\overline{k})^{*}T_{f^{*}E/P_{X}})=\deg((\sigma^{*}T_{E/P})\otimes_{k}\overline{k})=\deg(\sigma^{*}T_{E/P})<0$. Thus $f^{*}E$ is not semistable either what is again a contradiction to our assumption.}\end{proof}

\begin{rem}
One should notice that the other direction is not correct in the two previous statements in general: If $E$ is semistable, then in general $f^{*}E$ may be not semistable.
In the case of vector bundles one knows from a result by Gieseker (see \cite[Lemma 1.1]{Gie}) that the inverse direction holds, if the morphism $f$ is finite and separable.
In the case of a $G$-torsor whose structure group $G$ is a connected reductive algebraic group $G$ of finite presentation over $k$ the inverse direction holds following \cite[Lemma 6.8]{BH}, if $f$ is again finite and separable, but only in that case that $G$ has a representation ``of low height''.
\end{rem}

\subsection{The degree of a torsor} 

We need still the notion of a degree of a torsor:

\begin{defn}[\protect{\cite[Definition 3.2]{HN}}]
Let $X$ be a smooth projective curve of genus $g$ over a field $k$, $G$ a connected reductive algebraic group over $k$ and $E$ a $G$-torsor on $X$. Then we call the homomorphism $d_{E}\colon\Hom(G,k^{*})\to\mathbb{Z}, \chi\mapsto\deg(E_{\chi})$ the degree of $E$. Here $E_{\chi}=E\times_{\spec\,k}\mathbb{G}_{m}/G$ is the line bundle associated to $E$ by the character $\chi$, where the action of $G$ on $E\times_{\spec\,k}\mathbb{G}_{m}$ is given by \mbox{$(e,f)g=(eg,\chi(g^{-1})f)$} for all $e\in E$, $f\in\mathbb{G}_{m}$ and $g\in G$.
\end{defn}

\begin{rem}\label{null}
Let $X$ still be a smooth projective curve over a field $k$ and $G$ a connected reductive algebraic group over $k$. Then it holds: \\ $(i)$ The degree $d_{G_{X}}$ of the trivial torsor $G_{X}$ on $X$ is the zero homomorphism.\\
$(ii)$ If $E$ is a $G$-torsor on $X$ and $f:Y\to X$ a finite morphism of smooth projective curves over $k$, then we have $d_{f^{*}E}=[K(Y):K(X)]d_{E}$.\\
$(iii)$ If $E$ is a $G$-torsor on $X$ and $\overline{k}$ an algebraic closure of $k$, it holds \mbox{$d_{E\otimes_{k}\overline{k}}=d_{E}$}. 
\end{rem}

Important for us is this notion of a degree of a torsor mainly because of the following result by Holla and Narasimhan:

\begin{prop}[\protect{\cite[Theorem 1.2]{HN}}]
Let $G$ be a connected reductive algebraic group and $X$ a smooth projective curve both over an algebraic closed field $k$ of arbitrary characteristic. Then the set of isomorphism classes of semistable $G$-torsors on $X$ with fixed degree is bounded.\\
Here we say that a set $\mathcal{S}$ of $G$-torsors is bounded if there exists a scheme $S$ of finite type over $k$ and a family of $G$-torsors which is parametrized by $S$ such that any element of $\mathcal{S}$ is on $X$ isomorphic to the $G$-torsor on $X$ that one obtains by restriction of the given family on a suitable closed point of $S$.
\end{prop}

As by a theorem of Lang (see \cite[Theorem 1.3]{La}) \mbox{$H^{1}(k,G_{\overline{k}})=0$} for any connected reductive algebraic group $G$ over a finite field $k$, where $\overline{k}$ is an algebraic closure of $k$, such that there are no nontrivial forms, it follows by noetherian descent:

\begin{cor}\label{isoklassen}
Let $G$ be a connected reductive algebraic group and $X$ a smooth projective curve each over a finite field $k$. Then the set of isomorphism classes of semistable $G$-torsors on $X$ with fixed degree is bounded.
\end{cor}

\begin{cor}\label{endliso}
Let $G$ be a connected reductive algebraic group and $X$ a smooth projective curve each over a finite field $k$. Then there are only finitely many isomorphism classes of semistable $G$-torsors on $X$ with fixed degree.
\end{cor}

\subsection{Relations between semistability and trivialisibility}

Similary to the case of vector bundles which was examined by Lange and Stuhler in \cite{LS} we find also in the case of $G$-torsors relations between semistability and trivialisibility of torsors. Therefore we define:

\begin{defn}
Let $G$ be a connected reductive algebraic group of finite presentation over a perfect field $k$ of characteristic $p>0$ and $X$ a smooth projective curve over $k$. Then we call a $G$-torsor $E$ on $X$ \emph{strongly semistable}, if $E$ is semistable and moreover for all $r\in\N$ the pullback $(F^{r})^{*}E$ under the $r$-th product of the absolute Frobenius $F$ on $X$ is again semistable.
\end{defn}

This definition allows us to show the following result which is analogous to the well known result of Lange and Stuhler in the case of vector bundles (\cite[Satz 1.9.]{LS}):

\begin{prop}\label{ls19} 
Let $X$ be a smooth projective curve of genus $g$ over a finite field $k$ and $G$ a connected reductive algebraic group of finite presentation over $k$.
Then $E$ is, considered as a $G$-torsor for the \'etale topology, trivialisable, if and only if $E$ is strongly semistable of degree zero.
\end{prop}  
\begin{proof}{ If $E$ can be trivialized by a finite covering \mbox{$f\colon Y\to X$}, then for every $r\in\N$ the pullback $(F^{r})^{*}E$ under the $r$-th power of the absolute Frobenius on $X$ is also trivialisable by $f$. This implies that $(F^{r})^{*}E$ is semistable by Lemma \ref{endlichsemist} because the trivial torsor $G_{Y}$ on $Y$ is semistable by Corollary \ref{torsortrivial}. Thus $E$ is strongly semistable. From Remark \ref{null} follows moreover that $d_{E}\equiv0$ holds.\\
If conversely $E$ is strongly semistable of degree zero, then \mbox{$(F^{r})^{*}(E)$} is semistable of degree zero for all $r\in\N$. As by Corollary \ref{endliso} the number of isomorphism classes of semistable $G$-torsors $\widetilde{E}$ on $X$ with fixed degree $d_{\widetilde{E}}$ is finite, there exists an $s\in\N$ such that we have isomorphisms of $G$-torsors \mbox{$(F^{r})^{*}(E)\cong(F^{s})^{*}(E)$} for all $r\geq s$. By Corollary \ref{potenz}, $E$ is thus trivialisable.}\end{proof}

\subsection{Application to the characterisation of the category $\mathscr{B}_{\mathfrak{X}_{\mathfrak{o}},D}(G)$}

Let $G$ in the following always be a connected reductive group scheme of finite presentation over $\mathfrak{o}$ and $k=\overline{\mathbb{F}_{p}}$ the residue field of $\overline{\mathbb{Z}}_{p}$. Then we can now state our first main result:

\begin{thm}\label{thm20}
Let $\mathfrak{X}$ be a smooth model over $\overline{\mathbb{Z}}_{p}$ of a smooth projective curve $X$ over $\overline{\mathbb{Q}}_{p}$ of genus different from zero. A $G$-torsor $P$ on $\mathfrak{X}$ lies in the category $\mathscr{B}_{\mathfrak{X}_{\mathfrak{o}}}(G)$, if and only if $P_{k}$ is strongly semistable of degree zero on the smooth projective curve $\mathfrak{X}_{k}$ over $k$.
\end{thm}
\iproof{\begin{itemize}
\item[(i)] If $P\in\mathscr{B}_{\mathfrak{X}_{\mathfrak{o}}}(G)$, there exists by Theorem \ref{thm16} a covering \mbox{$(\pi\colon\mathcal{Y}\to\mathfrak{X})\in S_{\mathfrak{X}}^{\text{ss}}$} such that $\pi_{k}^{*}P_{k}$ is trivial. As $\pi_{k}$ is finite, $P_{k}$ is thus a semistable $G$-torsor on $\mathfrak{X}_{k}$ following Corollary \ref{torsortrivial} and Lemma \ref{endlichsemist}. As \mbox{$\pi\circ F_{\mathfrak{X}_{k}}=F_{\mathcal{Y}_{k}}\circ\pi$} holds, it follows that $(F_{\mathfrak{X}_{k}}^{r})^{*}P_{k}$ is also semistable for all \mbox{$r\in\N$} such that therefore $P_{k}$ is strongly semistable. As the degree of the trivial torsor is the zero morphism and \mbox{$d_{\pi_{k}^{*}P_{k}}=[K(\mathcal{Y}_{k}):K(\mathfrak{X}_{k})] d_{P_{k}}$} following Remark \ref{null}$(ii)$, the degree of $P_{k}$ is the zero morphism.

\item[(ii)] Let now conversely $P_{k}$ be strongly semistable of degree zero. As under the given assumptions every $G$-torsor is also a torsor with respect to the \'etale topology, we may consider $P$ and $P_{k}$  in the following always as torsors for the \'etale topology. From noetherian descent follows that the family $(X,\mathfrak{X},G,P_{k})$ descends to a family $(X_{K},\mathfrak{X}_{\mathfrak{o}_{K}},G_{0},P_{0})$ over a finite extension $K$ of $\overline{\mathbb{Q}}_{p}$ with residue field $\kappa\cong\mathbb{F}_{q},q=p^{r}$. Here $X_{K}$ is a smooth projective curve over $K$, $\mathfrak{X}_{\mathfrak{o}_{K}}$ a smooth model of $X_{K}$ over $\mathfrak{o}_{K}$, $G_{0}$ a connected reductive group of finite presentation over a normal finitely generated $\mathfrak{o}_{K}$-algebra $A$ and $P_{0}$ an \'etale $G$-torsor on the special fibre $\mathfrak{X}_{0}$ of $\mathfrak{X}_{\mathfrak{o}_{K}}$. One checks easily that $P_{0}$ is again strongly semistable of degree zero. By Proposition \ref{ls19} there exists thus a finite surjective morphism \mbox{$\varphi\colon\mathcal{Y}_{0}\to\mathfrak{X}_{0}$} such that $\varphi^{*}P_{0}$ is trivial. Following the construction of the morphism $\varphi$ in the proof of Proposition \ref{ls19} we may choose $\varphi$ as the composition of a finite \'etale morphism $\phi_{0}\colon\mathcal{Y}_{0}\to\mathfrak{X}_{0}$, where $\mathcal{Y}_{0}$ is a smooth projective curve over $\kappa$, and a Frobenius-power $F^{s}$ for a suitable $s\geq 0$.
From this point on, the remaining proof continues verbaly as those of \cite[Theorem 20]{De-We1}, if we only replace $\mathscr{E}_{k}$ by $P_{k}$ and the citation of \cite[Theorem 16]{De-We1} by those of Theorem \ref{thm16}.\endproof\end{itemize}

For the general case we define similary as in the case of vector bundles in \cite{De-We1}:

\begin{defn}
Let $R$ be a valuation ring with quotient field $Q$ and residue field $k$. We consider a model $\mathfrak{X}$ over $R$ of a smooth projective curve $X$ over $Q$ and a $G$-torsor $E$ on $\mathfrak{X}$ where $G$ is a connected reductive group scheme of finite presentation over $R$. Then we say that $E$ has \emph{strongly semistable reduction of degree zero} if the pullback of $E_{k}$ to the normalisation $\widetilde{C}$ of every irreducible component $C$ (equipped with the reduced structure) of $\mathfrak{X}_{k}$ is strongly semistable of degree zero. Note that every $\widetilde{C}$ is a smooth projective curve over $k$. 
\end{defn}

Moreover, we need the following generalisation of \cite[Theorem 18]{De-We1}:

\begin{prop}\label{thm18}
Let $G$ be a connected reduktive group scheme of finite presentation over $\mathbb{F}_{q}$, $X$ a pure onedimensional proper scheme over $\mathbb{F}_{q}$ and $E$ a $G$-torsor on $X$. Then the following are equivalent:
\begin{itemize}
\item[(i)] The pullback of $E$ to the normalisation of every irreducible component of $X$ is strongly semistable of degree zero.
\item[(ii)] There exists a finite surjective morphism $\varphi\colon Y\to X$, where $Y$ is a pure onedimensional proper scheme over $\mathbb{F}_{q}$ such that $\varphi^{*}E$ is a trivial $G$-torsor.
\item[(iii)] The same as in $(ii)$, but with $\varphi$ as a composition $\varphi\colon Y\xrightarrow{F^{s}}Y\xrightarrow{\pi}X$ for some $s\geq 0$, where $\pi$ is finite, \'etale and surjective and $F=Fr_{q}=Fr_{p}^{r}$ is the $q=p^{r}$-linear Frobenius on $Y$. 
\end{itemize}
\end{prop}
\iproof{\begin{itemize}{
\item[(ii)$\Rightarrow$(i)] Every irreducible component $C$ of $X$ is finitely dominated by an irreducible component $D$ of $Y$. It follows that the pullback of $E$ to $\widetilde{C}$ is trivialized by a finite morphism $\widetilde{D}\to\widetilde{C}$. As we may check semistability following Lemma \ref{endlichsemist} after pullback along a finite covering and the absolute Frobenius is functorial, this implies $(i)$.
\item[(i)$\Rightarrow$(iii)] By Corollary \ref{endliso}, there exist only finitely many  isomorphism classes of semi\-stable $G$-torsors with fixed degree on each of the finitely many irreducible components of $X$. As a first step we deduce from this that there exist only finitely many isomorphism classes of \'etale $G$-torsors on $X$ which all pullbacks to the normalisations of the irreducible components of $X$ are semistable of degree zero.

Let us first assume that $X$ is reduced. Then the normalisation morphism \mbox{$\pi\colon\widetilde{X}=\coprod\widetilde{C}_{v}\to X$} is finite where $X=\cup_{v}C_{v}$ is the partition of $X$ in its irreducible components. Moreover, the natural morphism $\alpha\colon G_{X}\to\pi_{*}G_{\widetilde{X}}$ is injective: The normalisation morphism $\pi$ is finite and surjective and thus a covering for the $h$-topology on $X$ (\cite[Definition 3.1.2]{V}), in particular therefore an universal topological epimorphism. As every scheme $U$ that is \'etale over $X$ is reduced, the morphism $\pi_{U}\colon U\times_{X}\widetilde{X}\to U$ is therefore following \cite[Lemma 3.2.1]{V} a categorial epimorphism in the category of schemes such that the map \[\alpha(U)\colon\Gamma(U,G_{U})\to\Gamma(U\times_{X}\widetilde{X},G_{U\times_{X}\widetilde{X}})\]  is injective.

As the support of the cocernel of $\alpha\colon G_{X}\to\pi_{*}G_{\widetilde{X}}$ consists at last only  of the singular points of $X$, this cocernel is a skyscraper sheaf of sets, following \cite[Exercise 2.2.9]{Liu} therefore of the form $\prod_{x\in X^{sing}}i_{x*}S_{x}$. We make now the assertion that every set $S_{x}$ is finite:

We show that we have the following exact sequence of \'etale sheaves on $X$:
\[0\to G_{X} \to \pi_{*}G_{\widetilde{X}} \to \bigoplus_{x\in X^{sing}}i_{x*}(\bigoplus_{y\mathop{|}x}G(\kappa(y))/G(\kappa(x))). \] Let $i\colon x\to X$ be a closed point. Then one has a cartesian diagram \begin{equation*}
        \xymatrix@=3em{%
         x\times_{X}\widetilde{X} \ar[r]^{\widetilde{i_{x}}} \ar[d]^{\pi_{x}} & \widetilde{X} \ar[d]^{\pi} \\
         x \ar@{^{(}->}[r]_{i_{x}} & X}
           \end{equation*} and there exists a canonical base change morphism \mbox{$i_{x}^{*}\pi_{*}G_{\widetilde{X}}\stackrel{\sim}{\longrightarrow} \pi_{x*}\widetilde{i_{x}}^{*}G_{\widetilde{X}}$} by \cite[Satz II.6.4.2]{Ta} which is an isomorphism.

           This implies $(\pi_{*}G_{\widetilde{X}})_{x}=(\pi_{x*}\widetilde{i_{x}}^{*}G_{\widetilde{X}})_{x}$ and by \cite[Corollary II.3.5]{M} \[(\pi_{*}G_{\widetilde{X}})_{x}=\bigoplus_{y\mathop{|}x}G_{\widetilde{X},y}.\] By \cite[Remark II.2.9d)]{M}, these stalks are \mbox{$G_{X,x}=G(\kappa(x))$} and \mbox{$G_{\widetilde{X},y}=G(\kappa(y))$}. As the stalk of $\bigoplus_{x\in X^{sing}}i_{x*}(\bigoplus_{y\mathop{|}x}G(\kappa(y))/G(\kappa(x)))$ in the point $x$ can be calculated to $\prod_{\widetilde{x}\mathop{|}x}(G(\kappa(\widetilde{x})/G(\kappa(x))$, the sequence \[0\to G_{X} \to \pi_{*}G_{\widetilde{X}} \to \bigoplus_{x\in X^{sing}}i_{x*}(\bigoplus_{y\mathop{|}x}G(\kappa(y))/G(\kappa(x)))\] is exact on any stalk and thus also itself exact. In particular follows therefore $S_{x}\subset\prod_{\widetilde{x}\mathop{|}x}(G(\kappa(\widetilde{x}))/G(\kappa(x)))$ for all $x\in X^{sing}$. From the finiteness of $\kappa(\widetilde{x})$ and the finiteness of the product it follows therefore that $S_{x}$ is finite for all $x\in X^{Sing}$.
           
As in the case of vector bundles it follows thus from \cite[Corollaire III.3.2.4]{Gi} that there are only finitely many isomorphism classes of $G$-torsors on $X$ which induce distinct chosen isomorphism classes of $G$-torsors on the curves $\widetilde{C}_{v}$. In particular, there exist thus only finitely many isomorphism classes of $G$-torsors on $X$ whose pullback to the normalisations of the irreducible components of $X$ is semistable of degree zero.
           
In the general case that $X$ is not necessarily reduced, it suffices therefore to show that the natural map $\varphi\colon H^{1}_{\acute{e}t}(X,G)\rightarrow H^{1}_{\acute{e}t}(X_{red},G)$ has finite fibres. By devissage, one sees that it is sufficient to show that for any ideal $\mathcal{J}\subset\mathcal{O}_{X}$ with $\mathcal{J}^{2}=0$ the map $i^{*}\colon H^{1}_{\acute{e}t}(X,G)\rightarrow H^{1}_{\acute{e}t}(X',G)$ has finite fibres where $i\colon X'\into X$ is the closed subscheme of $X$ defined by $\mathcal{J}$. 

Similar to the proof of Theorem \ref{thm16}, we have by \cite[Lemme VII.1.3.5]{Gi} the following exact sequence of sheaves on $X_{\acute{e}t}$ \begin{equation*}
        \xymatrix@=1em{%
         1 \ar[r] & i_{*}(Lie(G_{X'}/X')\otimes_{\mathcal{O}_{X}}\mathcal{J}) \ar[rrrr]^-{h} & & & & G \ar[rrr]^-{adj}& & & i_{*}G_{X'} \ar[r] & 1,}
           \end{equation*}
where all maps are defined as in the proof of Theorem \ref{thm16}. 

From this exact sequence one deduces the exact cohomology sequence \[H^{1}_{\acute{e}t}(X,i_{*}(Lie(G_{X'}/X')\otimes_{\mathcal{O}_{X'}}\mathcal{J}))\xrightarrow{h} H^{1}_{\acute{e}t}(X,G_{X})\to H^{1}_{\acute{e}t}(X,i_{*}G_{X'}).\]
By \cite[Corollaire VIII.5.6]{SGA 4}, it is moreover $R^{1}i_{*}G_{X'}=0$ because $i$ is a closed immersion, in particular integral. By \cite[Proposition V.3.1.3]{Gi}, $H^{1}_{\acute{e}t}(X,i_{*}G_{X'})$ is therefore isomorphic to $H^{1}_{\acute{e}t}(X',G_{X'})$ such that therefore the sequence  \[H^{1}_{\acute{e}t}(X,i_{*}(Lie(G_{X'}/X')\otimes_{\mathcal{O}_{X'}}\mathcal{J}))\xrightarrow{h} H^{1}_{\acute{e}t}(X,G_{X})\xrightarrow{i^{*}} H^{1}_{\acute{e}t}(X',G_{X'})\] is exact.

As above, \mbox{$R^{1}i_{*}(Lie(G_{X'}/X')\otimes_{\mathcal{O}_{X'}}\mathcal{J})=0$} implies by \cite[Proposition V.3.1.3]{Gi} that we have an isomorphism \begin{equation*}
 \xymatrix@=3em{%
H^{1}_{\acute{e}t}(X',Lie(G_{X'}/X')\otimes_{\mathcal{O}_{X'}}\mathcal{J})\ar[r]^{j}_{\sim} & H^{1}_{\acute{e}t}(X,i_{*}(Lie(G_{X'}/X')\otimes_{\mathcal{O}_{X'}}\mathcal{J}))}\end{equation*} such that the sequence \[H^{1}_{\acute{e}t}(X',Lie(G_{X'}/X')\otimes_{\mathcal{O}_{X'}}\mathcal{J})\stackrel{h\circ j}{\longrightarrow}H^{1}_{\acute{e}t}(X,G_{X})\xrightarrow{i^{*}} H^{1}_{\acute{e}t}(X',G_{X'})\] is exact.

The canonical $\mathcal{O}_{X'}$-module structure of $Lie(G_{X'}/X')$ induces moreover the following canonical exact sequence of sheaves on $X'_{\acute{e}t}$:
\[0\longrightarrow\ker(\alpha)\longrightarrow\mathcal{J}\stackrel{\alpha}{\longrightarrow}Lie(G_{X'}/X')\otimes_{\mathcal{O}_{X'}}\mathcal{J}\longrightarrow0\] where $\alpha$ is the natural morphism. This gives us a surjection \begin{equation*}\xymatrix@=2em{%
H^{1}_{\acute{e}t}(X',\mathcal{J})\ar@{>>}[rrr]^-{\alpha} & & & H^{1}_{\acute{e}t}(X',Lie(G_{X'}/X')\otimes_{\mathcal{O}_{X'}}\mathcal{J}),}\end{equation*} because $X'$ is at last onedimensional and $H^{2}_{\acute{e}t}(X',\mathscr{F})\cong H^{2}(X',\mathscr{F}_{Zar})$ for any \'etale sheaf $\mathscr{F}$ on $X'$ by \cite[Proposition III.3.3]{M} where $\mathscr{F}_{Zar}$ is the restriction to the Zariski topology.
Therefore one obtains the following non-abelian cohomology sequence: $H^{1}_{\acute{e}t}(X',\mathcal{J})\xrightarrow{h\circ j\circ\alpha} H^{1}_{\acute{e}t}(X,G_{X})\xrightarrow{i^{*}} H^{1}_{\acute{e}t}(X',G_{X'})$.
As $H^{1}_{\acute{e}t}(X',\mathcal{J})$ is a finite dimensional $\mathbb{F}_{q}$-vector space and therefore finite, it follows that $\varphi$ has finite fibres.

Let now $E$ be a $G$-torsor as given in $(i)$. Then the pullbacks to $\widetilde{C}_{v}$ of all $G$-torsors $(F_{X}^{n})^{*}E$ on $X$ are each semistable of degree zero. As we have shown that there exist only finitely many isomorphism classes of such $G$-torsors on $X$, there exists an $s\geq 0$ such that $(F_{X}^{t})^{*}E\cong(F^{s}_{X})^{*}E$ holds for all $t\geq s$. For the $G$-torsor $E':=(F_{X}^{s})^{*}E$ holds therefore $(F^{r}_{X})^{*}E'=E'$ for all $r=t-s\geq1$. If we consider now $E$ and $E'$ as torsors with respect to the \'etale topology, then there exists by Corollary \ref{pullback2} a finite \'etale and surjective morphism $\pi\colon Y\to X$ such that $\pi^{*}E'$ is a trivial $G$-torsor on $Y$. In other words, \mbox{$(\pi\circ F_{Y}^{s})^{*}E=(F^{s}_{X}\circ\pi)^{*}E=\pi^{*}(F^{s}_{X})^{*}E$} is a trivial \'etale $G$-torsor on $Y$. As $Y$ is a pure onedimensional proper $\mathbb{F}_{q}$-scheme, we have therefore shown $(iii)$.}\endproof\end{itemize}}

As a generalisation of \cite[Theorem 17]{De-We1} we are thus able to prove:

\begin{thm}\label{thm17}
Let $\mathfrak{X}$ be a model over $\overline{\mathbb{Z}}_{p}$ of a smooth projective curve $X$ over $\overline{\mathbb{Q}}_{p}$ and $G$ a connected reductive group scheme of finite presentation over $\mathfrak{o}$.
A $G$-torsor $P$ on $\mathfrak{X}_{\mathfrak{o}}$ is contained in the category $\mathscr{B}_{\mathfrak{X}_{\mathfrak{o}},D}(G)$ for a suitable divisor $D$ if and only if $P$ has strongly semistable reduction of degree zero. If this condition is fullfilled, there exist two divisors $D$ and $\widetilde{D}$ on $X$ with disjoint support such that $P$ lies as in the category $\mathscr{B}_{\mathfrak{X}_{\mathfrak{o}},D}(G)$ as in the category $\mathscr{B}_{\mathfrak{X}_{\mathfrak{o}},\widetilde{D}}(G)$.
\end{thm}
\begin{proof}{ If a $G$-torsor $P$ on $\mathfrak{X}_{\mathfrak{o}}$ is contained in the category $\mathscr{B}_{\mathfrak{X}_{\mathfrak{o}},D}(G)$ for a suitable divisor $D$ on $X$, then there exists by Theorem \ref{thm16} a covering $\pi\colon\mathcal{Y}\to\mathfrak{X}$ in $\Sh_{\mathfrak{X},D}^{\text{good}}$ such that $\pi_{k}^{*}P_{k}$ is trivial. Let $\mathfrak{X}_{k}=\cup_{v}C_{v}$ be the partition of $\mathfrak{X}_{k}$ in its irreducible components. Then every $C_{v}$ is finitely dominated by a irreducible component $\widetilde{C}_{v}$ of $\mathcal{Y}_{k}$. Thus it follows from Corollary \ref{torsortrivial} and Lemma \ref{endlichsemist} that the pullbacks of $P_{k}$ to all $\widetilde{C}_{v}$ are each strongly semistable of degree zero.

Let us now conversely assume that the $G$-torsor $P$ on $\mathfrak{X}_{\mathfrak{o}}$ has strongly semistable reduction of degree zero and is considered in the following always as an \'etale $G$-torsor. As in the proofs of Theorems \ref{thm16} and \ref{thm20}, there exists by noetherian descent a finite extension $K$ of $\mathbb{Q}_{p}$ with ring of integers $\mathfrak{o}_{K}$ and residue field $\kappa\cong\mathbb{F}_{q}$ such that the family $(X,\mathfrak{X},C_{v},G,P_{k})$ descends to a family $(X_{K},\mathfrak{X}_{\mathfrak{o}_{K}},C_{v0},G_{0},P_{0})$ with analogous properties. In particular, $P_{0}$ is an \'etale $G_{0}$-torsor on the special fibre $\mathfrak{X}_{0}=\mathfrak{X}_{\mathfrak{o}_{K}}\otimes\kappa$ whose pullback to the normalisations $\widetilde{C}_{v0}$ of the irreducible components $C_{v0}$ of $\mathfrak{X}_{0}$ is strongly semistable of degree zero.
By Proposition \ref{thm18}, we find therefore a finite \'etale morphism $\widetilde{\pi_{0}}\colon\widetilde{\mathcal{Y}_{0}}\to\mathfrak{X}_{0}$ such that for the composition $\widetilde{\varphi_{0}}\colon\widetilde{\mathcal{Y}_{0}} \xrightarrow{F^{s}}\widetilde{\mathcal{Y}_{0}}\xrightarrow{\widetilde{\pi_{0}}}\mathfrak{X}_{0}$ the pullback $\widetilde{\varphi_{0}}^{*}P_{0}$  is trivial. Note that we may replace here $s$ by any greater number $s'$ and then $F$ by any arbitrary power of $F$. As in \cite{De-We1}, \cite[Expos\'e IX, Th\'eor\`eme 1.10]{SGA1} allows us to lift the morphism $\widetilde{\pi_{0}}$ to a finite \'etale morphism $\widetilde{\pi_{\mathfrak{o}_{K}}}\colon\widetilde{\mathcal{Y}_{\mathfrak{o}_{K}}}\to\mathfrak{X}_{0}$ with $\widetilde{\pi_{0}}$ as special fibre. If we replace $K$ as in \cite{De-We1} by a suitable finite extension $K_{1}$ of $K$, we can dominate $\widetilde{\pi_{\mathfrak{o}_{K_{1}}}}$ by an object $\pi_{\mathfrak{o}_{K_{1}}}\colon\mathcal{Y}_{\mathfrak{o}_{K_{1}}}\to\mathfrak{X}_{\mathfrak{o}_{K_{1}}}$ of the category $\Sh_{\mathfrak{X}_{\mathfrak{o}_{K_{1}}},\emptyset}$, where we may suppose as in \cite{De-We1} that $\mathcal{Y}_{\mathfrak{o}_{K_{1}}}$ is not only even semistable and moreover also regular. If we write again $K$ instead of $K_{1}$, $P_{0}$ instead of $P_{1}$, $\mathfrak{o}_{K}$ instead of $\mathfrak{o}_{K_{1}}$, $\kappa$ instead of $\kappa_{1}$ etc. to simplify our notation, it follows that the pullback of $P_{0}$ under the composition $\varphi\colon\mathcal{Y}_{0}\xrightarrow{F^{s}}\mathcal{Y}_{0}\xrightarrow{\pi_{0}}\mathfrak{X}_{0}$ is a trivial $G_{0}$-torsor.

From this point on, the remaining proof carries over verbatim from the proof of \cite[Theorem 17]{De-We1}, if we only replace $\mathcal{E}_{\kappa}$ by $P_{k}$, $\mathcal{E}$ by $P$ and $\mathcal{E}_{0}$ by $P_{0}$ and finally the citation of \cite[Theorem 16]{De-We1} by those of Theorem \ref{thm16}. }\end{proof}

\section{Principal bundles on $X_{\mathbb{C}_{p}}$ and parallel transport}
\subsection{Potentially strongly semistable reduction and parallel transport}

Let $X$ still be a smooth projective curve over $\overline{\mathbb{Q}}_{p}$, $\mathfrak{X}$ a model of $X$ over $\overline{\mathbb{Z}}_{p}$ and $G$ a connected reductive group scheme of finite presentation over $\mathfrak{o}$. In the aim to extend the results of the previous section to the case of principal bundles over $X_{\mathbb{C}_{p}}$, we need the following definition which is analogous to those in the case of vector bundles (see \cite[Chapter 0]{De-We1}):

\begin{defn}
We say that a $G$-torsor $E$ on $X_{\mathbb{C}_{p}}$ has \emph{strongly semistable reduction of degree zero}, if he has the following property:\\
$E$ extends to a $G$-torsor $\widetilde{E}$ on $\mathfrak{X}_{\mathfrak{o}}$ for a suitable model $\mathfrak{X}$ of $X$ over $\overline{\mathbb{Z}}_{p}$ and the pullback of the special fibre $\widetilde{E}_{k}$ of $\widetilde{E}$ to the normalisation of every irreducible component of $\mathfrak{X}_{k}$ is strongly semistable of degree zero.\\
Moreover, one says that $E$ has \emph{potentially strongly semistable reduction of degree zero}, if there exists a finite \'etale morphism $\alpha\colon Y\to X$ of smooth projective curves over $\overline{\mathbb{Q}}_{p}$ such that $\alpha^{*}_{\mathbb{C}_{p}}E$ has strongly semistable reduction of degree zero.   
\end{defn}

Let $\mathcal{B}_{X_{\mathbb{C}_{p}}}^{s}(G)$ now be the category of $G$-torsors on $X_{\mathbb{C}_{p}}$ with strongly semistable reduction of degree zero and $\mathcal{B}_{X_{\mathbb{C}_{p}}}^{ps}(G)$ the category of those with potentially strongly semistable reduction of degree zero. Then we are now able to show the two main results of this article:

\begin{thm}\label{s}
It holds $\mathcal{B}^{s}_{X_{\mathbb{C}_{p}}}(G)=\bigcup_{D}\mathcal{B}_{X_{\mathbb{C}_{p}}, D}(G)$. Every $G$-torsor in $\mathcal{B}^{s}_{X_{\mathbb{C}_{p}}}(G)$ is contained as in $\mathcal{B}_{X_{\mathbb{C}_{p}}, D}(G)$ as in $\mathcal{B}_{X_{\mathbb{C}_{p}}, \widetilde{D}}(G)$ for suitable divisors $D$ and $\widetilde{D}$ with disjoint support.

There exists an unique continous functor \mbox{$\rho_{E}\colon\pi_{1}(X)\to\mathcal{P}(G(\mathbb{C}_{p}))$} for every $G$-torsor $E\in\mathcal{B}^{s}_{X_{\mathbb{C}_{p}}, D}(G)$ such that $\rho_{E}(x)=E_{x}$ for all \mbox{$x\in X(\mathbb{C}_{p})$} and $\rho_{E}$ is compatible with the functors $\rho_{E}\colon\pi_{1}(X-D)\to\mathcal{P}(G(\mathbb{C}_{p}))$ which were constructed for fixed divisors $D$ with $E\in\mathcal{B}_{X_{\mathbb{C}_{p}},D}(G)$ in the second section.

In conclusion one obtains therefore a functor $\rho\colon\mathcal{B}_{X_{\mathbb{C}_{p}}}^{s}(G)\to\Rep_{\pi_{1}(X)}(G(\mathbb{C}_{p}))$ which is functorial with respect to morphisms of smooth projective curves over $\overline{\mathbb{Q}}_{p}$, morphisms of connected reductive group schemes of finite presentation over $\mathfrak{o}$ and $\mathbb{Q}_{p}$-automorphisms of $\overline{\mathbb{Q}}_{p}$ and which is compatible with the analogous functors in the case of vector bundles through the functors defined in Section 2.4.

For any point $x_{0}\in X(\mathbb{C}_{p})$ is the ``fibre functor in $x_{0}$'' \[\mathcal{B}_{X_{\mathbb{C}_{p}}}^{s}(G)\to \mathcal{P}(G(\mathbb{C}_{p})),\,E\mapsto E_{x_{0}},\,f\mapsto f_{x_{0}}\] faithful.
\end{thm}
\begin{proof}{ If $E\in\mathcal{B}_{X_{\mathbb{C}_{p}},D}(G)$ for some divisor $D$, then $E$ extends by the definition of the category $\mathcal{B}_{X_{\mathbb{C}_{p}},D}(G)$ and by Theorem \ref{thm17} to a $G$-torsor $\widetilde{E}$ on $\mathfrak{X}_{\mathfrak{o}}$ for a suitable model $\mathfrak{X}$ on $X$ which has strongly semistable reduction of degree zero. By definition, it is therefore $E\in\mathcal{B}^{s}_{X_{\mathbb{C}_{p}}}(G)$.

If conversely $E'\in\mathcal{B}^{s}_{X_{\mathbb{C}_{p}}}(G)$ holds, then it follows either from Theorem \ref{thm17} that $E'\in\mathcal{B}_{X_{\mathbb{C}_{p}},D}(G)$ for a suitable divisor $D$ and that there are more precisely divisors $D$ and $\widetilde{D}$ with disjoint support such that $E$ is contained as in $\mathcal{B}_{X_{\mathbb{C}_{p}}, D}(G)$ as in $\mathcal{B}_{X_{\mathbb{C}_{p}}, \widetilde{D}}(G)$. This implies $\mathcal{B}^{s}_{X_{\mathbb{C}_{p}}}(G)=\bigcup_{D}\mathcal{B}_{X_{\mathbb{C}_{p}}, D}(G)$.

If $E$ is a $G$-torsor in $\mathcal{B}_{X_{\mathbb{C}_{p}}}^{s}(G)$, there exist, for all divisors $D$ with \mbox{$E\in\mathcal{B}_{X_{\mathbb{C}_{p}},D}(G)$}, by Section 2.2 continous functors $\rho_{E,D}\colon\mathcal{B}_{X_{\mathbb{C}_{p}},D}(G)\to\Rep_{\pi_{1}(X-D)}(G(\mathbb{C}_{p}))$ with $\rho_{E}(x)=E_{x}$ for all $x\in (X-D)(\mathbb{C}_{p})$. By Proposition \ref{seifert}, there exists thus an unique functor $\rho_ {E}\colon\pi_{1}(X)\to\mathcal{P}(G(\mathbb{C}_{p}))$ which induces the functors $\rho_{E,D}$ on \mbox{$\pi_{1}(X-D)$}.

It is clear that we obtain a functor $\rho\colon\mathcal{B}_{X_{\mathbb{C}_{p}}}^{s}(G)\to\Rep_{\pi_{1}(X)}(G(\mathbb{C}_{p}))$ by the application $E\mapsto\rho_{E}$ and all functoriality properties follow directly from the analogous statements in Section 2.4.

Let us fix a point $x_{0}\in X(\mathbb{C}_{p})$. We have to show the injectivity of the natural map \mbox{$\Hom_{\mathcal{B}_{X_{\mathbb{C}_{p}}}^{s}(G)}(E,E')\rightarrow\Hom_{\mathcal{P}(G(\mathbb{C}_{p}))}(E_{x_{0}},E'_{x_{0}})$} for all \mbox{$G$-torsors $E,\,E'\in\mathcal{B}_{X_{\mathbb{C}_{p}}}^{s}(G)$}: Let \mbox{$f,g\colon E\rightarrow E'$} be two morphisms of $G$-torsors on $X_{\mathbb{C}_{p}}$ such that \mbox{$f_{x_{0}}=g_{x_{0}}$} holds. Let moreover \mbox{$x\in X(\mathbb{C}_{p})$} be another $\mathbb{C}_{p}$-valued point of $X$ and $\gamma$ a path from $x_{0}$ to $x$. Then there exist two parallel transports $\rho_{E}(\gamma)\colon E_{x_{0}}\rightarrow E_{x}$ resp.  \mbox{$\rho_{E'}(\gamma)\colon E'_{x_{0}}\rightarrow E'_{x}$} such that we obtain by the functoriality of the construction of the parallel transport a commutative diagram \begin{equation*}
        \xymatrix@=3em{%
         E_{x_{0}} \ar[r]^{\rho_{E}(\gamma)} \ar[d] & E_{x}\ar[d] \\ 
       E'_{x_{0}} \ar[r]^{\rho_{E'}(\gamma)} & E'_{x} }
           \end{equation*} as in the case that one chooses as vertical arrows $f_{x_{0}}$ and $f_{x}$ as either in the case that one chooses $g_{x_{0}}$ and $g_{x}$. By the commutativity of the diagram and the fact that all horizontal maps are isomorphisms, it follows $f_{x}=g_{x}$. As $x\in X(\mathbb{C}_{p})$ was chosen arbitrarily, $f_{x}=g_{x}$ holds thus for all $\mathbb{C}_{p}$-valued points $x$ of $X$ and in particular for all $\mathbb{C}_{p}$-valued points of $E$. Therefore the maps $\widetilde{f},\widetilde{g}\colon E(\mathbb{C}_{p})\rightarrow E'(\mathbb{C}_{p})$ induced by $f$ resp. $g$ coincidence. As $G_{X_{\mathbb{C}_{p}}}$ is a connected reductive algebraic group of finite presentation over $X_{\mathbb{C}_{p}}$, we know that $E$ and $E'$ are affine, smooth, reductive and of finite presentation over $X_{\mathbb{C}_{p}}$ and therefore by \cite[Proposition 5.3.4]{EGAII} projective over $\mathbb{C}_{p}$. Thus $E$ and $E'$ come by \cite[p. 50]{EH} from varieties $\mathcal{E}$ resp. $\mathcal{E'}$ over $\mathbb{C}_{p}$ in the sense of classical algebraic geometry \`a la Weil. As a morphism $\mathcal{E}\to\mathcal{E'}$ is uniquement determined by a morphism $\mathcal{E}(\mathbb{C}_{p})\to\mathcal{E'}(\mathbb{C}_{p})$, then $f$ and $g$ coincidence on the set of closed points of $E$ which is an open dense subset of $E$ by \cite[Exercise II.3.14]{Ha}. Following \cite[Lemma 7.2.2.1]{EGA1} $f$ and $g$ therefore conincidence on the whole scheme $E$ such that the fibre functor in $x_{0}$ is faithful.}\end{proof}

We can now prove the theorem mentioned in the introduction, but recall first \cite[Proposition 31]{De-We1} that we need for the proof:

\begin{prop}[\protect{\cite[Proposition 31]{De-We1}}]\label{31}
Let \mbox{$\alpha\colon Y\to X$} be a Galois covering between varieties over $\overline{\mathbb{Q}}_{p}$. A (continous) functor $W\colon\pi_{1}(Y)\to\mathcal{C}$ into a (topological) category $\mathcal{C}$ factors as $W=V\circ\alpha_{*}$ for some (continous) functor $V\colon\pi_{1}(X)\to\mathcal{C}$ if and only if we have $W\circ\sigma_{*}=W$ for all $\sigma\in\Gal(Y/X)$. If $\alpha$ is only finite and \'etale, but not necessarily Galois, the relation $W=V\circ\alpha_{*}$ already determines $V$ uniquely.
\end{prop}

\begin{thm}\label{main}
Let $E$ be a $G$-torsor in $\mathcal{B}_{X_{\mathbb{C}_{p}}}^{ps}(G)$, i.\,e. a $G$-torsor on $X_{\mathbb{C}_{p}}$ with potentially strongly semistable reduction of degree zero. There exist functorial isomorphisms of ``parallel transport'' along \'etale paths of fibres of $E_{\mathbb{C}_{p}}$ on $X_{\mathbb{C}_{p}}$. In particular, there exists for any $G$-torsor $E\in\mathcal{B}^{s}_{X_{\mathbb{C}_{p}}, D}(G)$ an unique continous functor $\rho_{E}\colon\pi_{1}(X)\to\mathcal{P}(G(\mathbb{C}_{p}))$ with $\rho_{E}(x)=E_{x}$ for all \mbox{$x\in X(\mathbb{C}_{p})$} such that $\rho_{E}$ is compatible with the above defined functors $\rho_{\alpha^{*}E}\colon\pi_{1}(Y)\to\mathcal{P}(G(\mathbb{C}_{p}))$, if $\alpha^{*}E\in\mathcal{B}_{Y_{\mathbb{C}_{p}}}^{s}(G)$ holds for some finite \'etale morphism \mbox{$\alpha\colon Y\to X$} between smooth projective curves over $\overline{\mathbb{Q}}_{p}$.

In this way, one gets a functor $\rho\colon\mathcal{B}_{X_{\mathbb{C}_{p}}}^{ps}(G)\to\Rep_{\pi_{1}(X)}(G(\mathbb{C}_{p}))$ which is functorial with respect to morphisms of smooth projective curves over $\overline{\mathbb{Q}}_{p}$, morphisms of connected reductive group schemes of finite presentation over $\mathfrak{o}$ and $\mathbb{Q}_{p}$-automorphisms of $\overline{\mathbb{Q}}_{p}$ and which is compatible via the functors described in Section 2.4. with the analogous functors in the case of vector bundles.

For every point $x_{0}\in X(\mathbb{C}_{p})$, the ``fibre functor in $x_{0}$'' \[\mathcal{B}_{X_{\mathbb{C}_{p}}}^{ps}\to \mathcal{P}(G(\mathbb{C}_{p})), E\mapsto E_{x_{0}}, f\mapsto f_{x_{0}}\] is faithful.
\end{thm}
\begin{proof}{ Let $E$ be a $G$-torsor on $X_{\mathbb{C}_{p}}$ with potentially strongly semistable reduction of degree zero. Then there exists a finite \'etale morphism $\alpha\colon Y\to X$ of smooth projective curves over $\overline{\mathbb{Q}}_{p}$ such that $E\in\mathcal{B}_{Y_{\mathbb{C}_{p}}}^{s}(G)$. Without any loss of generality we may assume that $\alpha$ is Galois. By Theorem \ref{s} it holds therefore \mbox{$\rho_{\alpha^{*}E}\circ\sigma_{*}=\rho_{\sigma^{*}(\alpha^{*}E)}=\rho_{\alpha^{*}E}$} for all \mbox{$\sigma\in\Gal(Y/X)$}. Following Proposition \ref{31}, there exists thus an unique continous functor $\rho(E)=\rho_{E}\colon\pi_{1}(X)\rightarrow\mathcal{P}(G(\mathbb{C}_{p}))$ such that $\rho_{\alpha{*}E}=\rho_{E}\circ\alpha_{*}$. In particular, we have $\rho_{E}(x)=E_{x}$ for all \mbox{$x\in X(\mathbb{C}_{p})$} and a morphism $\rho_{E}(\gamma)=\rho_{\alpha^{*}E}(\gamma')\colon E_{x_{1}}=(\alpha^{*}E)_{y_{1}}\rightarrow(\alpha^{*}E)_{y_{2}}=E_{x_{2}}$ for any \'etale path $\gamma$ from $x_{1}$ to $x_{2}$ in $X$. Here \mbox{$y_{1}\in Y(\mathbb{C}_{p})$} is lying over $x_{1}$ and $\gamma'$ is the unique path in $Y$ from $y_{1}$ to the point $y_{2}$ over $x_{2}$ satisfying \mbox{$\alpha_{*}\gamma'=\gamma$}. For a morphism $f\colon E\to E'$ of $G$-torsors in $\mathcal{B}^{ps}_{X_{\mathbb{C}_{p}}}(G)$, the morphism \mbox{$\rho(f)=\rho_{f}\colon\rho_{E}\to\rho_{E'}$} is defined by the family of maps $f_{x}\colon E_{x}\to E'_{x}$ for all $x\in X(\mathbb{C}_{p})$.

We assert now that this construction defines a well-defined functor \[\rho\colon\mathcal{B}^{ps}_{X_{\mathbb{C}_{p}}}(G)\longrightarrow\Rep_{\pi_{1}(X)}(G(\mathbb{C}_{p}))\] which extends the previous functors $\rho$ on $\mathcal{B}_{X_{\mathbb{C}_{p}}}^{s}(G)$. From the functoriality properties in Theorem \ref{s} and the uniqueness statement in Proposition \ref{31} follows that the definition of $\rho_{E}$ is independent of the choice of $\alpha$.

Next we show that for any morphism $f\colon E\rightarrow E'$ in $\mathcal{B}_{X_{\mathbb{C}_{p}}}^{ps}(G)$ the family of maps $f_{x}\colon E_{x}\rightarrow E'_{x}$ defines a morphism in $\Rep_{\pi_{1}(X)}(G(\mathbb{C}_{p}))$. Without any loss of generality, we may assume that both $\alpha^{*}E$ and $\alpha^{*}E'$ are contained in $\mathcal{B}_{Y_{\mathbb{C}_{p}}}^{s}(G)$. Then $\rho_{\alpha^{*}f}$, i.\,e. the family of maps $(\alpha^{*}f)_{y}\colon(\alpha^{*}E)_{y}\rightarrow(\alpha^{*}E')_{y}$, defines a morphism in $\Rep_{\pi_{1}(Y)}(G(\mathbb{C}_{p}))$. By the definition of $\rho_{E}(\gamma)$ for an \'etale path in $X$, the diagram \begin{equation*}
        \xymatrix@=3em{%
         E_{x_{1}} \ar[r]^{f_{x_{1}}} \ar[d]^{\rho_{E}(\gamma)} & E'_{x_{1}}\ar[d]^{\rho_{E'}(\gamma)} \\ 
       E_{x_{2}} \ar[r]^{f_{x_{2}}} & E'_{x_{2}} }
           \end{equation*} commutes for any \'etale path $\gamma$ in $X$ such that therefore the family $(f_{x})_{x\in X(\mathbb{C}_{p})}$ defines a morphism in $\Rep_{\pi_{1}(X)}(G(\mathbb{C}_{p}))$.
           
Thus $\rho$ is well-defined and it is obvious by the construction that $\rho$ is a functor and extends the functor given on $\mathcal{B}_{X_{\mathbb{C}_{p}}}^{s}(G)$. All the claimed functorialities of $\rho_{E}$ follow from the already known functorialities of $\rho_{\alpha^{*}E}$ and the uniqueness statement in Proposition \ref{31}. The fact that the fibre functor $\mathcal{B}_{X_{\mathbb{C}_{p}}}^{ps}(G)\to\mathcal{P}(G(\mathbb{C}_{p}))$ is faithful, is shown in the same way as in Theorem \ref{s}. }\end{proof}

\subsection{Open questions}

\begin{itemize}
\item How big is the category $\mathcal{B}_{X_{\mathbb{C}_{p}}}^{ps}(G)$? Is in particular every semistable $G$-torsor of degree zero on $X_{\mathbb{C}_{p}}$ contained in this category? The analogous question in the case of vector bundles, i.\,e. whether any semistable vector bundle of degree zero on $X_{\mathbb{C}_{p}}$ lies in $\mathcal{B}^{ps}_{X_{\mathbb{C}_{p}}}$, is until now solved by a postive answer only for curves $X$ of genus $g=0$ or $g=1$ and is still open for curves $X$ of higher genus. 
\item What is the essential image of the functor $\rho\colon\mathcal{B}_{X_{\mathbb{C}_{p}}}^{ps}(G)\to\Rep_{\pi_{1}(X)}(G(\mathbb{C}_{p}))$? This is still an open question in the case of vector bundles either.
\item Is the fibre functor in a fixed point $x_{0}\in X(\mathbb{C}_{p})$, i.\,e. $\mathcal{B}_{X_{\mathbb{C}_{p}}}^{ps}(G)\to\mathcal{P}(G(\mathbb{C}_{p}))$, not only faithful, but moreover fully faithful?
\end{itemize}

\textit{Westf\"alische Wilhelms-Universit\"at M\"unster,
Mathematisches Institut,
Einsteinstrasse 62,
D-48149 M\"unster, Germany\\
current adress: Universit\"at Ulm, Institut f\"ur Reine Mathematik, Helmholtzstrasse 18, D-89081 Ulm, Germany, urs.hackstein@uni-ulm.de}

\begin{thebibliography}{999}
\bibitem [Be]{Be}K. A. Behrend: \textit{Semi-stability of reductive group schemes over curves.} Math. Ann. {\bf 301}, no. 2 (1995), p. 281--305
\bibitem[BH]{BH}I. Biswas and Y. Holla: \textit{Harder-Narasimham reduction of a principal bundle.} Nagoya Math. J. {\bf 174} (2004), p. 201--223
\bibitem[De-We1]{De-We1} C. Deninger and A. Werner: \textit{Vector bundles on $p$-adic curves and parallel transport.} Ann. Scient. \'Ec. Norm. Sup. {\bf 38} (2005), p. 535--597
\bibitem [De-We2] {De-We2}C. Deninger and A. Werner: \textit{On Tannaka duality for vector bundles on $p$-adic curves.} In: J. Nagel and C. Peters (eds): Algebraic cycles and Motives. London Mathematical Society  Lecture Note Series {\bf 344}(2006)
\bibitem [De-We3] {De-We3} C. Deninger and A. Werner: \textit{Principal bundles on $p$-adic curves and parallel transport.} Private Notes (2005) 
\bibitem[EGAI]{EGA1}A. Grothendieck and  J. Dieudonn\'e: \textit{\'El\'ements de G\'eometrie Alg\'ebrique I.} Publ. Math. IHES {\bf 4} (1960)
\bibitem[EGAII]{EGAII}A. Grothendieck and J. Dieudonn\'e: \textsl{\'El\'ements de G\'eometrie Alg\'ebrique II.} Publ. Math. IHES {\bf 8} (1961)
\bibitem[EH]{EH}D. Eisenbud and J. Harris: \textit{The Geometry of schemes.} Graduate Texts in Mathematics {\bf 197}, Springer, Berlin, Heidelberg, New York (2000) 
\bibitem[FK]{FK}E. Freitag and R. Kiehl: \textit{\'Etale cohomology and the Weil conjecture.} Ergebnisse der Mathematik und ihrer Grenzgebiete, 3. Folge, Band {\bf 13}, Springer, Berlin, Heidelberg, New York (1988)
\bibitem[Gi]{Gi}J. Giraud: \textit{Cohomologie non ab\'elienne.} Springer, Berlin, Heidelberg, New York (1971)
\bibitem [Gie]{Gie}D. Gieseker: \textit{On a theorem of Bogomolov on Chern classes of stable bundles.} Am. J. Math. {\bf 101} (1979), p. 79--85
\bibitem[Gr]{Gr} A. Grothendieck: \textit{Le groupe de Brauer III: Exemples et complements.} In: J. Giraud, A. Grothendieck, S.L. Kleiman, M. Raynaud and J. Tate: Dix expos\'es sur la cohomologie des sch\'emas. North-Holland Publishing Company, Amsterdam, Masson\&Cie, Paris (1968), p. 88--188
\bibitem[H]{H} U. Hackstein: \textit{Prinzipalb\"undel auf $p$-adischen Kurven und Paralleltransport.} PhD-thesis, WWU M\"unster (2006), published online on www.miami-uni-muenster.de and arXiv:math.AG/0701315
\bibitem[Ha]{Ha}R. Hartshorne: \textit{Algebraic Geometry.} Springer, Berlin, Heidelberg, New York (1977)
\bibitem[He]{He}G. Herz: \textit{On representations attached to semistable vector bundles on Mumford curves.} PhD-thesis WWU M\"unster (2005)
\bibitem[HN]{HN}Y. Holla and M. Narasimhan: \textit{A generalisation of Nagata's theorem on ruled surfaces.} Compositio Mathematica {\bf 127} (2001), p. 321--332
\bibitem [La]{La}S. Lang: \textit{Algebraic groups over finite fields.} Am. J. Math. {\bf 78} (1956), p. 555--563
\bibitem[Las]{Las}Y. Laszlo: \textit{A non-trivial family of bundles fixed by the square of Frobenius.} C. R. Acad. Sci. Paris, t. {\bf 333}, S\'erie I (2001), p.651--656
\bibitem[LS]{LS}H. Lange and U. Stuhler:\textit{ Vektorb\"undel auf Kurven und Darstellungen der algebraischen Fundamentalgruppe.} Math. Z. {\bf 156} (1977), 73--83
\bibitem[Liu1]{Liu}Q. Liu: \textit{Algebraic Geometry and Arithmetic Curves.} Oxford University Press (2002)
\bibitem[Liu2]{Liu2}Q. Liu: \textit{Stable reduction of finite covers of curves.} arXiv:math.AG/0412075v2 (5 Jan. 2005)
\bibitem[M]{M}J. Milne: \textit{\'Etale cohomology.} Princeton University Press, Princeton, New Jersey (1980)
\bibitem[Mf]{Mf} D. Mumford: \textit{Abelian varieties.} Oxford University Press (1970)
\bibitem[Mu]{Mu} J.P. Murre: \textit{Lectures on an introduction to Grothendieck's theory of the fundamental group.} Tata institute of fundamental research, Bombay (1967)
\bibitem[Na-Se]{Na-Se}M.S. Narasimhan and C.S. Seshadri: \textit{Stable and unitary vector bundles on a compact Riemann surface.} Ann. Math. {\bf 82} (1965), 540--567
\bibitem[R]{R}A. Ramanathan: \textit{Stable principal bundles on a compact Riemann surface.} Math. Ann. {\bf 213} (1975), p. 129--152
\bibitem [Se]{Se} J.-P. Serre: \textit{Espaces fibr\'es alg\'ebriques.} In: S\'eminaire C. Chevalley ENS, 2e ann\'ee 1958: Anneaux de Chow et applications (1958), p. 1--37
\bibitem[SGA1]{SGA1}A. Grothendieck et. al.: \textit{S\'eminaire de g\'eometrie alg\'ebrique du Bois-Marie 1960/61. Rev\^etements \'etales et groupe fondamental. }LNM {\bf 224}, Springer, Berlin, Heidelberg, New York (1971)
\bibitem[SGA3]{SGA3}M. Demazure and A. Grothendieck: \textit{S\'eminaire de g\'eometrie alg\'ebrique du Bois-Marie. Sch\'emas en groupes.} LNM {\bf 151}, Springer, Berlin, Heidelberg, New York (1970)
\bibitem [SGA4]{SGA 4}M. Artin, A. Grothendieck and J.L. Verdier (eds.):\textit{ S\'eminaire de G\'eom\'etrie alg\'ebrique du Bois-Marie 1963-1964: Theorie des topos et cohomologie \'etale des sch\'emas.} LNM {\bf 269}, {\bf 270}, Springer, Berlin, Heidelberg, New York (1972)
\bibitem [T]{T} J. Tong: \textit{Application d'Albanese pour les courbes et contractions.} Pr\'e\-pub\-li\-ca\-tion, Orsay (2006)
\bibitem [Ta]{Ta}G. Tamme: \textit{Einf\"uhrung in die \'etale Kohomologie.} Der Regensburger Trichter {\bf 17}, Universit\"at Regensburg (1979)
\bibitem [V]{V} V. Voevodsky: \textit{Homology of schemes.} Selecta Mathematica, New Series, Vol. {\bf 2}, No.1 (1996), p. 111--153
\bibitem[W]{W}A. Weil: \textit{G\'en\'eralisation des fonctions ab\'eliennes.} J. de Math. P. et App. (IX) {\bf 17} (1938) 47--87
\bibitem [We]{We}A. Werner: \textit{Vector bundles on curves over $\mathbb{C}_{p}$. }in: L. Weng and I. Nakamura (eds.): Arithmetic geometry and number theory. Series in number theory and its applications, Vol. {\bf 1}, World Scientific (2006), p. 47--62
\end{thebibliography}
\end{document}